\documentclass[11pt]{article}
\usepackage{amsmath,amssymb,amsbsy,amsfonts,latexsym,amsthm,mathrsfs}
\usepackage{color}
\usepackage{graphicx}
\usepackage{float}
\newtheorem{theorem}{Theorem}[section]
\newtheorem{proposition}[theorem]{Proposition}
\newtheorem{definition}[theorem]{Definition}

\newtheorem{lemma}{Lemma}[section]
\newtheorem{remark}{Remark}[section]

\input pictex.sty
\begin{document}
\title{Asymptotic results for certain first-passage times and areas
of renewal processes\thanks{The authors acknowledge the support
of GNAMPA-INdAM and of MIUR Excellence Department Project awarded
to the Department of Mathematics, University of Rome Tor Vergata
(CUP E83C18000100006).}}
\author{Claudio Macci\thanks{Address: Dipartimento di Matematica,
Universit\`a di Roma Tor Vergata, Via della Ricerca Scientifica,
I-00133 Rome, Italy. e-mail: \texttt{macci@mat.uniroma2.it}}\and
Barbara Pacchiarotti\thanks{Address: Dipartimento di Matematica,
Universit\`a di Roma Tor Vergata, Via della Ricerca Scientifica,
I-00133 Rome, Italy. e-mail: \texttt{pacchiar@mat.uniroma2.it}}}
\date{}
\maketitle
\begin{abstract}
\noindent We consider the process $\{x-N(t):t\geq 0\}$, where $x\in\mathbb{R}_+$
and $\{N(t):t\geq 0\}$ is a renewal process with light-tailed distributed holding times. We are 
interested in the joint distribution of $(\tau(x),A(x))$ where $\tau(x)$ is the first-passage
time of $\{x-N(t):t\geq 0\}$ to reach zero or a negative value, and 
$A(x):=\int_0^{\tau(x)}(x-N(t))dt$ is the corresponding first-passage (positive) area swept out 
by the process $\{x-N(t):t\geq 0\}$. We remark that we can define the sequence 
$\{(\tau(n),A(n)):n\geq 1\}$ by referring to the concept of integrated random walk. Our aim is 
to prove asymptotic results as $x\to\infty$ in the fashion of large (and moderate) deviations.
\end{abstract}

\section{Introduction}
Let $\{N(t):t\geq 0\}$ be the renewal process defined by
\begin{equation}\label{eq:renewal-process}
N(t):=\sum_{n\geq 1}1_{T_1+\cdots+T_n\leq t},
\end{equation}
where $\{T_n:n\geq 1\}$ be i.i.d. positive random variables; then, for $x\in\mathbb{R}_+$,
let $\tau(x)$ be the first-passage time of $\{x-N(t):t\geq 0\}$ to reach zero or a negative value, and let 
$A(x)$ is the corresponding first-passage (positive) area swept out by the process $\{x-N(t):t\geq 0\}$, 
i.e.
$$A(x):=\int_0^{\tau(x)}(x-N(t))dt.$$

Here we generalize the presentation in \cite{AbundoFuria}, where $\{N(t):t\geq 0\}$ is
a Poisson process. However, according to the terminology in some other references in the literature,
we can refer to the concept of \emph{integrated random walk}, at least when $x$ is integer.
In fact we can consider the random walk $\{S_n:n\geq 1\}$ defined by $S_n:=\sum_{j=1}^nT_j$, and the
bivariate sequence $\{(S_n,S_1+\cdots+S_n):n\geq 1\}$ coincides with the sequence $\{(\tau(n),A(n)):n\geq 1\}$
presented above. Among the references with asymptotic results for integrated random walks here we
recall \cite{BorovkovBoxmaPalmowski} and \cite{DenisovPerfilevWatchel} for the heavy-tailed case, and
\cite{PerfilevWatchel} for the light-tailed case. 
Here, for completeness, we also recall \cite{Book} and \cite{Saulis} which provide exact large deviation
probabilities for more general weighted random walks.

Throughout this paper the random variables $\{T_n:n\geq 1\}$ are assumed to be
light-tailed distributed; this allows us to apply the G\"artner Ellis Theorem (see e.g. 
Theorem \ref{th:GE}recalled in this paper), and we can obtain the large deviation principle for
$\left\{\left(\frac{\tau(x)}{x},\frac{A(x)}{x^2}\right):x>0\right\}$ (as $x\to\infty$) under suitable
hypotheses. The asymptotic bounds provided by the large deviation principle allow us to 
estimate the exponential decay rate of probabilities of rare events by computing the infimum of the rate
function over suitable sets (see Remark \ref{rem:convergences} for more details).

As far as we know this work is the first attempt to study the asymptotic behavior of the
bivariate family $\left\{\left(\frac{\tau(x)}{x},\frac{A(x)}{x^2}\right):x>0\right\}$; in a successive work
one could try to obtain \emph{exact} asymptotic results as in the ones in \cite{PerfilevWatchel} for the
marginal distributions only.

In this paper we also study moderate deviations. More precisely we mean a class of large deviation principles
for families of random variables depending on the choice of certain scalings in a suitable class; all these
large deviation principles (whose speed function depends on the scaling) are governed by the same quadratic
rate function vanishing at the origin. In some sense this class of large deviation principles fills the gap
between two asymptotic regimes, i.e. the convergence of $\left(\frac{\tau(x)}{x},\frac{A(x)}{x^2}\right)$
to a constant as $x\to\infty$, and an asymptotic normality result (this will be explained in Remark
\ref{rem:MD-typical-feature}).

In some literature the results on the distribution of first-passage areas
(possibly in terms of the joint distribution with first-passage times) concern Markov processes,
and in particular some L\'evy processes; see e.g. the jump-diffusion processes in \cite{Abundo}
and the drifted Brownian motion \cite{AbundoDelvescovo}. This approach allows to consider
suitable differential-difference equations (in terms of the generator of the Markov process) for
the Laplace transform of $(\tau(x),A(x))$ which can be solved. However, if
$\{N(t):t\geq 0\}$ is a renewal process, it is easy to check that the random variable $A(x)$
can be expressed in terms of a suitable linear combinations of the holding times. So, in this
case, we can easily manage the joint distribution of $(\tau(x),A(x))$ even if
$\{N(t):t\geq 0\}$ is not a Markov process (it is well-known that a renewal process is
Markovian if and only if the holding times are exponentially distributed).

The study of first-passage areas is often motivated by potential applications. For
instance first-passage areas can model the evolution of certain random systems described
by diffusive continuous processes (some examples in the literature concern solar physics
studies, non-oriented animal movement patterns, and DNA breathing dynamics) or by their
superpositions with jump processes; among the references in the literature we recall
\cite{KearneyMajumdarMartin}, \cite{KearneyPyeMartin}, \cite{Knight} and
\cite{PermanWellner} for processes without jumps, and \cite{KouWang} for a process with
jumps (more precisely in that reference there is a compound Poisson process with double
exponentially distributed jumps). Some further applications of first-passage areas
concern the framework of default-at-maturity models in finance (see e.g. \cite{Bates}
and \cite{Jorion}). Finally another application in queueing theory is described in
\cite{Abundo} where $\tau(x)$ is interpreted as the busy period (that is the time until
the queue is first empty) and $A(x)$ represents the cumulative waiting time experienced
by all the \lq\lq customers\rq\rq\ during a busy period.

For completeness we also cite \cite{Arendarczyk-et-al} and \cite{BlanchetMandjes} in which,
for suitable functions $T(u)$ and for some classes of processes $\{Q(t):t\geq 0\}$ having
interest in queueing theory, the authors study the asymptotic behavior of
$P\left(\int_0^{T(u)}Q(r)dr>u\right)$, as $u\to\infty$, in the fashion of large deviations.

We conclude with the outline of the paper. We start with some preliminaries in Section
\ref{sec:preliminaries}. We study large and moderate deviations in Sections
\ref{sec:LD} and \ref{sec:MD}, respectively.
In Section \ref{sec:ex} we present some details for the Poisson process case, i.e. for the 
case in which the random variables $\{T_n:n\geq 1\}$ are exponentially distributed. Finally
Section \ref{sec:conclusions} is devoted to some conclusions.

\section{Preliminaries}\label{sec:preliminaries}
In this section we present some preliminaries 
on large and moderate deviations (Section \ref{sub:preliminaries-ldmd}), on the joint 
distribution of $(\tau(x),A(x))$ (Section \ref{sub:joint-distribution-tau-A}) and on 
a suitable function $\Lambda$ which plays a crucial role in our results (Section 
\ref{sub:function-Lambda}).

\subsection{On large and moderate deviations}\label{sub:preliminaries-ldmd}
We start with the definition of large deviation principle (LDP from now on). 
In view of what follows we consider the LDP (as $x\to\infty$) for a family of $\mathbb{R}^h$-valued 
random variables $\{Z_x:x>0\}$ defined on the same probability space $(\Omega,\mathcal{F},P)$.

A lower semi-continuous function $I:\mathbb{R}^h\to[0,\infty]$ is called
rate function, and it is said to be good if all its level sets
$\{\{z\in\mathbb{R}^h:I(z)\leq\eta\}:\eta\geq 0\}$ are compact.
Then $\{Z_x:x>0\}$ satisfies the LDP with speed $v_x\to\infty$ and rate function $I$ if
$$\limsup_{x\to\infty}\frac{1}{v_x}\log P(Z_x\in C)\leq-\inf_{z\in C}I(z)\ \mbox{for all closed sets}\ C$$
and
$$\liminf_{x\to\infty}\frac{1}{v_x}\log P(Z_x\in O)\geq-\inf_{z\in O}I(z)\ \mbox{for all open sets}\ O.$$

We talk about moderate deviations when we have a class of LDPs for 
families of centered (or asymptotically centered) random variables which depends
on some scaling factors and, moreover, all these LDPs (whose speed functions
depend on the scaling factors) are governed by the same quadratic rate
function vanishing at zero. We can also say that, as usually happens, this
class of LDPs fills the gap between a convergence to a constant and an
asymptotic normality result; this will be illustrated in Remark
\ref{rem:MD-typical-feature} (a version of this remark can also be adapted to
the random variables in Remark \ref{rem:MD-same-result-with-different-centering-ee}).

The main large deviation tool used in this paper is the G\"artner Ellis
Theorem (see e.g. Theorem 2.3.6 in \cite{DemboZeitouni}), and here we
recall its statement.

\begin{theorem}\label{th:GE}
Assume that, for all $\alpha\in\mathbb{R}^h$, there exists
$$f(\alpha):=\lim_{x\to\infty}\frac{1}{v_x}\log\mathbb{E}\left[e^{v_x\langle\alpha,Z_x\rangle}\right]$$
as an extended real number (here $\langle\cdot,\cdot\rangle$ is the inner
product in $\mathbb{R}^h$); moreover assume that the origin $\alpha=0$
belongs to the interior of the set
$$\mathcal{D}(f):=\{\alpha\in\mathbb{R}^h:f(\alpha)<\infty\}.$$
Furthermore let $f^*$ be the function defined by
$$f^*(z):=\sup_{\alpha\in\mathbb{R}^h}\{\langle\alpha,z\rangle-f(\alpha)\}.$$
Then:
(a) for all closed sets $C$
$$\limsup_{x\to\infty}\frac{1}{v_x}\log P(Z_x\in C)\leq-\inf_{z\in C}f^*(z);$$
(b) for all open sets $O$
$$\liminf_{x\to\infty}\frac{1}{v_x}\log P(Z_x\in O)\geq-\inf_{z\in O\cap\mathcal{E}}f^*(z),$$
where $\mathcal{E}$ is the set of exposed points of $f^*$ (roughly speaking
is the set where $f^*$ is finite and strictly convex; see e.g. Definition
2.3.3 in \cite{DemboZeitouni});\\
(c) if $f$ is essentially smooth and lower semi-continuous, then the LDP holds.
\end{theorem}

For completeness we also recall that $f$ is essentially smooth (see e.g. Definition
2.3.5 in \cite{DemboZeitouni}) if the interior of $\mathcal{D}(f)$ is non-empty, if
$f$ is differentiable throughout the interior of $\mathcal{D}(f)$, and if $f$ is steep
(namely if $\|\nabla f(\alpha)\|\to\infty$ as $\alpha$ converges to a boundary
point $\alpha^{(0)}$ of $\mathcal{D}(f)$).

\subsection{On the joint distribution of $(\tau(x),A(x))$}\label{sub:joint-distribution-tau-A}
Let $\{N(t):t\geq 0\}$ be the renewal process defined by eq. \eqref{eq:renewal-process}
above, where the holding times $\{T_n:n\geq 1\}$ are i.i.d. positive random variables;
then their (common) moment generating function is
$$\mathbb{E}\left[e^{\alpha T_n}\right]=e^{\varphi(\alpha)}\ \mbox{for all}\ \alpha\in\mathbb{R},$$
for a suitable increasing function $\varphi$. Obviously we have $\varphi(0)=0$
and $\varphi(\alpha)<\infty$ for all $\alpha\leq 0$. Throughout this paper we always
assume that the origin $\alpha=0$ belongs to the interior of the set
$$\mathcal{D}(\varphi):=\{\alpha\in\mathbb{R}:\varphi(\alpha)<\infty\};$$
so we assume to have one of the two following cases:
\begin{itemize}
\item $\mathcal{D}(\varphi)=\mathbb{R}$;
\item $\mathcal{D}(\varphi)=(-\infty,\bar{\alpha})$ or $\mathcal{D}(\varphi)=(-\infty,\bar{\alpha}]$ for some $\bar{\alpha}>0$.
\end{itemize}
For instance, if $\{T_n:n\geq 1\}$ are exponentially distributed, i.e.
$\{N(t):t\geq 0\}$ is a Poisson process, we have
\begin{equation}\label{eq:varphi-exponential-case}
\varphi(\alpha):=\left\{\begin{array}{ll}
\log\frac{\lambda}{\lambda-\alpha}&\ \mbox{if}\ \alpha<\lambda\\
\infty&\ \mbox{otherwise}
\end{array}\right.\ \mbox{for some}\ \lambda>0;
\end{equation}
thus $\mathcal{D}(\varphi)=(-\infty,\bar{\alpha})$ and $\bar{\alpha}=\lambda$.

Furthermore let $\{X(t):t\geq 0\}$ be the process defined by
$$X(t):=x-N(t),\ \mbox{for some}\ x>0.$$
Obviously we have $X(0)=x$. In this paper we are interested in the joint distribution
of $(\tau(x),A(x))$ where
$$\tau(x):=\inf\{t\geq 0:X(t)\leq 0\}$$
is the first-passage time to reach zero or a negative value, and
$$A(x):=\int_0^{\tau(x)}X(t)dt$$
is the corresponding first-passage area. In particular we need to refer to the
moment generating functions and we have two cases; so we present two lemmas. We remark that 
we recover eqs. (3.17)-(3.18) in \cite{AbundoFuria} stated for the case of Poisson process,
i.e. for the case in which the function $\varphi$ is defined by \eqref{eq:varphi-exponential-case},
for some $\lambda>0$. From now on we use the notation $[x]:=\max\{k\in\mathbb{Z}:k\leq x\}$.

\begin{lemma}\label{lem:mgf-AbundoFuria-(3.17)-gen}
If $x$ is integer we have
$$\mathbb{E}\left[e^{\alpha_1\tau(x)+\alpha_2A(x)}\right]=\left\{\begin{array}{ll}
\prod_{k=1}^xe^{\varphi(\alpha_1+\alpha_2k)}&\ \mbox{if}\ (\alpha_1,\alpha_2)\in\mathcal{D}_x^{(1)}\\
\infty&\ \mbox{otherwise},
\end{array}\right.$$
where
$$\mathcal{D}_x^{(1)}:=\{(\alpha_1,\alpha_2)\in\mathbb{R}^2:\alpha_2\geq 0,\alpha_1+\alpha_2x\in\mathcal{D}(\varphi)\}
\cup\{(\alpha_1,\alpha_2)\in\mathbb{R}^2:\alpha_2<0,\alpha_1+\alpha_2\in\mathcal{D}(\varphi)\}.$$
\end{lemma}
\begin{proof}
If $x$ is integer we have
$$(\tau(x),A(x))=\left(\sum_{k=1}^xT_k,\sum_{k=1}^x(x-k+1)T_k\right).$$
Then, since $(\tau(x),A(x))$ and $\left(\sum_{k=1}^xT_k,\sum_{k=1}^xkT_k\right)$ are
identically distributed by the hypotheses, the joint moment generating function is
$$\mathbb{E}\left[e^{\alpha_1\tau(x)+\alpha_2A(x)}\right]=\left\{\begin{array}{ll}
\prod_{k=1}^xe^{\varphi(\alpha_1+\alpha_2k)}
&\ \mbox{if $\alpha_1+\alpha_2k\in\mathcal{D}(\varphi)$ for all $k\in\{1,\ldots,x\}$}\\
\infty&\ \mbox{otherwise};
\end{array}\right.$$
thus the desired equality holds by taking into account the definition of the set
$\mathcal{D}_x^{(1)}$ in the statement.
\end{proof}

\begin{lemma}\label{lem:mgf-AbundoFuria-(3.18)-gen}
If $x$ is not integer we have
$$\mathbb{E}\left[e^{\alpha_1\tau(x)+\alpha_2A(x)}\right]=\left\{\begin{array}{ll}
\prod_{k=0}^{[x]}e^{\varphi(\alpha_1+\alpha_2(x-k))}&\ \mbox{if}\ (\alpha_1,\alpha_2)\in\mathcal{D}_x^{(2)}\\
\infty&\ \mbox{otherwise},
\end{array}\right.$$
where
$$\mathcal{D}_x^{(2)}:=\{(\alpha_1,\alpha_2)\in\mathbb{R}^2:\alpha_2\geq 0,\alpha_1+\alpha_2x\in\mathcal{D}(\varphi)\}
\cup\{(\alpha_1,\alpha_2)\in\mathbb{R}^2:\alpha_2<0,\alpha_1+\alpha_2(x-[x])\in\mathcal{D}(\varphi)\}.$$
\end{lemma}
\begin{proof}
If $x$ is not integer we have
$$(\tau(x),A(x))=\left(\sum_{j=1}^{[x]+1}T_j,\sum_{j=1}^{[x]+1}(x-j+1)T_k\right).$$
Then, by setting $k=j-1$, the joint moment generating function is
$$\mathbb{E}\left[e^{\alpha_1\tau(x)+\alpha_2A(x)}\right]=\left\{\begin{array}{ll}
\prod_{k=0}^{[x]}e^{\varphi(\alpha_1+\alpha_2(x-k))}
&\ \mbox{if $\alpha_1+\alpha_2(x-k)\in\mathcal{D}(\varphi)$ for all $k\in\{0,\ldots,[x]\}$}\\
\infty&\ \mbox{otherwise};
\end{array}\right.$$
thus the desired equality holds by taking into account the definition of the set
$\mathcal{D}_x^{(2)}$ in the statement.
\end{proof}

We remark that we have $\mathcal{D}_x^{(1)}=\mathbb{R}^2$ in Lemma \ref{lem:mgf-AbundoFuria-(3.17)-gen},
and $\mathcal{D}_x^{(2)}=\mathbb{R}^2$ in Lemma \ref{lem:mgf-AbundoFuria-(3.18)-gen}, if and only if
$\mathcal{D}(\varphi)=\mathbb{R}$.

\subsection{A suitable function $\Lambda$}\label{sub:function-Lambda}
Here we discuss some properties of a suitable function $\Lambda$ which plays a crucial role in our results.
We start with its definition.

\begin{definition}\label{def:function-Lambda}
We have the following two cases.
\begin{itemize}
\item If $\mathcal{D}(\varphi)=\mathbb{R}$, then
$$\Lambda(\alpha_1,\alpha_2):=\int_0^1\varphi(\alpha_1+\alpha_2y)dy.$$
\item If $\mathcal{D}(\varphi)=(-\infty,\bar{\alpha})$ or $\mathcal{D}(\varphi)=(-\infty,\bar{\alpha}]$ for some $\bar{\alpha}>0$, then
$$\Lambda(\alpha_1,\alpha_2)=\left\{\begin{array}{ll}
\int_0^1\varphi(\alpha_1+\alpha_2y)dy&\ \mbox{if}\ (\alpha_1,\alpha_2)\in\mathcal{D}\\
\infty&\ \mbox{otherwise},
\end{array}\right.$$
where
$$\mathcal{D}:=\{(\alpha_1,\alpha_2)\in\mathbb{R}^2:\alpha_2\geq 0,\alpha_1+\alpha_2\in\mathcal{D}(\varphi)\}
\cup\{(\alpha_1,\alpha_2)\in\mathbb{R}^2:\alpha_2<0,\alpha_1\leq\bar{\alpha}\}.$$
\end{itemize}
Moreover, in both cases, it is easy to check the following formulas:
$$\int_0^1\varphi(\alpha_1+\alpha_2y)dy=\left\{\begin{array}{ll}
\varphi(\alpha_1)&\ \mbox{if}\ \alpha_2=0\\
\frac{1}{\alpha_2}\int_{\alpha_1}^{\alpha_1+\alpha_2}\varphi(y)dy&\ \mbox{if}\ \alpha_2\neq 0.
\end{array}\right.$$
\end{definition}

Now a brief discussion on the set
$$\mathcal{D}(\Lambda):=\{(\alpha_1,\alpha_2)\in\mathbb{R}^2:\Lambda(\alpha_1,\alpha_2)<\infty\}.$$
Obviously we have $\mathcal{D}(\Lambda)=\mathbb{R}^2$ if $\mathcal{D}(\varphi)=\mathbb{R}$.
On the other hand, if $\mathcal{D}(\varphi)=(-\infty,\bar{\alpha})$ or $\mathcal{D}(\varphi)=(-\infty,\bar{\alpha}]$
for some $\bar{\alpha}>0$, then $\mathcal{D}(\Lambda)\subset\mathcal{D}$ and we can have
several different situations. In particular we have the following three cases and, for
each one, we also discuss the lower semi-continuity of $\Lambda$.
\begin{itemize}
\item If $\mathcal{D}(\varphi)=(-\infty,\bar{\alpha}]$ for some $\bar{\alpha}>0$, then 
$\mathcal{D}(\Lambda)=\mathcal{D}$. In this case the function $\Lambda$ is
lower semi-continuous. An example is the inverse Gaussian distribution, where 
$T_1$ has density
$$f_{T_1}(x)=\frac{1}{\sqrt{2\pi}x^{3/2}}\exp\left(-\frac{1}{2x}(\mu x-1)^2\right)\quad \mbox{for some}\ \mu>0;$$
then
$$\varphi(\alpha):=\left\{\begin{array}{ll}
\mu-\sqrt{\mu^2-2\alpha}&\ \mbox{if}\ \alpha\leq\frac{\mu^2}{2}\\
\infty&\ \mbox{otherwise}.
\end{array}\right.$$
\item If $\mathcal{D}(\varphi)=(-\infty,\bar{\alpha})$ for some $\bar{\alpha}>0$
and $\varphi$ is integrable in a left neighborhood of $\alpha=\bar{\alpha}$,
as happens for the exponential distribution, i.e. $\varphi$ is as in eq.
\eqref{eq:varphi-exponential-case}, then $\mathcal{D}(\Lambda)=\mathcal{D}$.
In this case the function $\Lambda$ is not lower semi-continuous; indeed, for
$(\alpha_1^{(0)},\alpha_2^{(0)})$ in the boundary of $\mathcal{D}$ with $\alpha_2^{(0)}>0$ (and
therefore $(\alpha_1^{(0)},\alpha_2^{(0)})=(\bar{\alpha}-\alpha_2^{(0)},\alpha_2^{(0)})$),
the condition
\begin{equation}\label{eq:lsc-condition}
\liminf_{(\alpha_1,\alpha_2)\to(\bar{\alpha}-\alpha_2^{(0)},\alpha_2^{(0)})}\Lambda(\alpha_1,\alpha_2)
\geq\Lambda(\bar{\alpha}-\alpha_2^{(0)},\alpha_2^{(0)})
\end{equation}
fails because the left hand side is finite and the right hand side is equal to
infinity.
\item If $\mathcal{D}(\varphi)=(-\infty,\bar{\alpha})$ for some $\bar{\alpha}>0$
and $\varphi$ is not integrable in a left neighborhood of $\alpha=\bar{\alpha}$,
then $\mathcal{D}(\Lambda)$ coincides with the interior of $\mathcal{D}$. In this case
the function $\Lambda$ is lower semi-continuous. An example is the non central chi squared 
distribution, where $T_1$ has density
$$f_{T_1}(x)=\frac{1}{2}e^{-(x+\lambda)/2}\left(\frac{x}{\lambda}\right)^{k/4-1/2}I_{k/2-1}(\sqrt{\lambda x})\quad \mbox{for some}\ \lambda,k>0,$$
and $I_\nu(y)=(y/2)^\nu\sum_{j=0}^\infty\frac{(y^2/4)^j}{j!\Gamma(\nu+j+1)}$ is the modified 
Bessel function of the first kind; then
$$\varphi(\alpha):=\left\{\begin{array}{ll}
\frac{\lambda\alpha}{1-2\alpha}-\frac{k}{2}\log(1-2\alpha)&\ \mbox{if}\ \alpha<\frac{1}{2}\\
\infty&\ \mbox{otherwise}.
\end{array}\right.$$
\end{itemize}

Now we take the partial derivatives of $\Lambda$ when $(\alpha_1,\alpha_2)$ belongs to
the interior of $\mathcal{D}(\Lambda)$. Then, after some computations, we get
$$\left(\frac{\partial\Lambda}{\partial\alpha_1}(\alpha_1,\alpha_2),\frac{\partial\Lambda}{\partial\alpha_2}(\alpha_1,\alpha_2)\right)
=\left\{\begin{array}{ll}
\left(\varphi^\prime(\alpha_1),\frac{1}{2}\varphi^\prime(\alpha_1)\right)&\ \mbox{if}\ \alpha_2=0\\
\left(\frac{\varphi(\alpha_1+\alpha_2)-\varphi(\alpha_1)}{\alpha_2},
\frac{\alpha_2\varphi(\alpha_1+\alpha_2)-\int_{\alpha_1}^{\alpha_1+\alpha_2}\varphi(y)dy}{\alpha_2^2}\right)&\ \mbox{if}\ \alpha_2\neq 0;
\end{array}\right.$$
so, in view of what follows, we recall that
\begin{equation}\label{eq:gradient-of-Lambda-at-the-origin}
\left(\frac{\partial\Lambda}{\partial\alpha_1}(0,0),\frac{\partial\Lambda}{\partial\alpha_2}(0,0)\right)
=\left(\varphi^\prime(0),\frac{1}{2}\varphi^\prime(0)\right).
\end{equation}
Moreover we have the following results.

\begin{lemma}\label{lem:differentiability-of-Lambda}
The function $\Lambda$ is differentiable throughout the interior of $\mathcal{D}(\Lambda)$.
\end{lemma}
\begin{proof}
The statement can be proved showing that the partial derivatives of $\Lambda$ are continuous.
The continuity of the partial derivatives can be easily checked. In particular, for the
continuity of $\frac{\partial\Lambda}{\partial\alpha_2}$ for $\alpha_2=0$, we have
$$\lim_{\alpha_2\to 0}\frac{\alpha_2\varphi(\alpha_1+\alpha_2)-\int_{\alpha_1}^{\alpha_1+\alpha_2}\varphi(y)dy}{\alpha_2^2}
=\frac{1}{2}\varphi^\prime(\alpha_1)$$
by applying the Hopital rule.
\end{proof}

\begin{lemma}\label{lem:steepness-for-Lambda}
Assume that $\mathcal{D}(\varphi)\neq\mathbb{R}$. Then:
\begin{enumerate}
\item If $\mathcal{D}(\varphi)=(-\infty,\bar{\alpha}]$ for some $\bar{\alpha}>0$,
then the function $\Lambda$ is not steep.
\item If $\mathcal{D}(\varphi)=(-\infty,\bar{\alpha})$ for some $\bar{\alpha}>0$,
then the function $\Lambda$ is steep.
\end{enumerate}
\end{lemma}
\begin{proof}
For the Statement 1 we can observe that, if we take $(\alpha_1,\alpha_2)$ in the
interior of $\mathcal{D}(\Lambda)$ converging to a boundary point
$(\alpha_1^{(0)},\alpha_2^{(0)})$ with $\alpha_2^{(0)}\neq 0$, then
$$\left(\frac{\partial\Lambda}{\partial\alpha_1}(\alpha_1,\alpha_2),\frac{\partial\Lambda}{\partial\alpha_2}(\alpha_1,\alpha_2)\right)\to
\left(\frac{\varphi(\alpha_1^{(0)}+\alpha_2^{(0)})-\varphi(\alpha_1^{(0)})}{\alpha_2^{(0)}},
\frac{\alpha_2^{(0)}\varphi(\alpha_1^{(0)}+\alpha_2^{(0)})-\int_{\alpha_1^{(0)}}^{\alpha_1^{(0)}+\alpha_2^{(0)}}\varphi(y)dy}{(\alpha_2^{(0)})^2}\right),$$
and therefore the partial derivatives do not diverge.

For the Statement 2 it is enough to check that
$\frac{\partial\Lambda}{\partial\alpha_1}(\alpha_1,\alpha_2)$ diverges. In fact
we have the following cases.
\begin{itemize}
\item If $\alpha_2^{(0)}>0$, then $\varphi(\alpha_1+\alpha_2)\to\varphi(\alpha_1^{(0)}+\alpha_2^{(0)})=\varphi(\bar{\alpha})=\infty$,
and $\varphi(\alpha_1^{(0)})<\infty$.
\item If $\alpha_2^{(0)}<0$, then $\varphi(\alpha_1)\to\varphi(\alpha_1^{(0)})=\varphi(\bar{\alpha})=\infty$,
and $\varphi(\alpha_1^{(0)}+\alpha_2^{(0)})<\infty$.
\item If $\alpha_2^{(0)}=0$, then $(\alpha_1,\alpha_2)\to (\alpha_1^{(0)},\alpha_2^{(0)})=(\bar{\alpha},0)$
and we have
$$\lim_{(\alpha_1,\alpha_2)\to(\bar{\alpha},0)}\frac{\partial\Lambda}{\partial\alpha_1}(\alpha_1,\alpha_2)
=\lim_{\alpha_1\to\bar{\alpha}}\varphi^\prime(\alpha_1)=\infty.$$
\end{itemize}
\end{proof}

We conclude with some further preliminaries. Let $C=(c_{ij})_{i,j\in\{1,2\}}$ be the
Hessian matrix of $\Lambda$ at the origin; then, after some computations, we get
\begin{equation}\label{eq:hessian-matrix-of-Lambda-at-the-origin}
C:=\left(\begin{array}{cc}
\frac{\partial^2\Lambda}{\partial\alpha_1^2}(0,0)&\frac{\partial^2\Lambda}{\partial\alpha_1\partial\alpha_2}(0,0)\\
\frac{\partial^2\Lambda}{\partial\alpha_2\partial\alpha_2}(0,0)&\frac{\partial^2\Lambda}{\partial\alpha_2^2}(0,0)
\end{array}\right)
=\varphi^{\prime\prime}(0)\left(\begin{array}{cc}
1&1/2\\
1/2&1/3
\end{array}\right).
\end{equation}
Moreover let $C^{-1}=(\hat{c}_{ij})_{i,j\in\{1,2\}}$ is the inverse of the
square matrix $C$ in eq. \eqref{eq:hessian-matrix-of-Lambda-at-the-origin},
and therefore
\begin{equation}\label{eq:inverse-of-hessian-matrix-of-Lambda-at-the-origin}
C^{-1}:=\frac{1}{\varphi^{\prime\prime}(0)}\left(\begin{array}{cc}
4&-6\\
-6&12
\end{array}\right).
\end{equation}

\section{Large deviations of $\left\{\left(\frac{\tau(x)}{x},\frac{A(x)}{x^2}\right):x>0\right\}$}\label{sec:LD}
We start with the first main result in this paper. 
The proof consists of an application of Theorem \ref{th:GE} with $h=2$, and the function $f$ is the function 
$\Lambda$ in Definition \ref{def:function-Lambda}.

\begin{proposition}\label{prop:LD}
Let $\Lambda$ be the function in Definition \ref{def:function-Lambda}, (which is finite in a
neighborhood of the origin $(\alpha_1,\alpha_2)=(0,0)$), and let $\Lambda^*$ be defined by
\begin{equation}\label{eq:def-Lambda-star}
\Lambda^*(z_1,z_2):=\sup_{(\alpha_1,\alpha_2)\in\mathbb{R}^2}\{\alpha_1z_1+\alpha_2z_2-\Lambda(\alpha_1,\alpha_2)\}.
\end{equation}
Then:
$$\limsup_{x\to\infty}\frac{1}{x}\log P\left(\left(\frac{\tau(x)}{x},\frac{A(x)}{x^2}\right)\in C\right)
\leq-\inf_{(z_1,z_2)\in C}\Lambda^*(z_1,z_2)\ \mbox{for all closed sets}\ C$$
and
$$\liminf_{x\to\infty}\frac{1}{x}\log P\left(\left(\frac{\tau(x)}{x},\frac{A(x)}{x^2}\right)\in O\right)
\geq-\inf_{(z_1,z_2)\in O\cap\mathcal{E}}\Lambda^*(z_1,z_2)\ \mbox{for all open sets}\ O,$$
where $\mathcal{E}$ is the set of exposed points of $\Lambda^*$.
\end{proposition}
\begin{proof}
We want to apply the G\"artner Ellis Theorem; so we have to show that
$$\lim_{x\to\infty}\frac{1}{x}\log\mathbb{E}\left[e^{x\left(\alpha_1\frac{\tau(x)}{x}+\alpha_2\frac{A(x)}{x^2}\right)}\right]
=\Lambda(\alpha_1,\alpha_2)\ (\mbox{for every}\ (\alpha_1,\alpha_2)\in\mathbb{R}^2),$$
where $\Lambda$ is the function in Definition \ref{def:function-Lambda}. In particular
we note that
$$\frac{1}{x}\log\mathbb{E}\left[e^{x\left(\alpha_1\frac{\tau(x)}{x}+\alpha_2\frac{A(x)}{x^2}\right)}\right]
=\frac{1}{x}\log\mathbb{E}\left[e^{\alpha_1\tau(x)+\frac{\alpha_2}{x}A(x)}\right].$$
In what follows we consider the cases $x$ integer and $x$ non-integer. In both cases, if
$\mathcal{D}(\varphi)=\mathbb{R}$, some parts below can be simplified; in fact the conditions
$(\alpha_1,\alpha_2/x)\in\mathcal{D}_x^{(1)}$ and $(\alpha_1,\alpha_2/x)\in\mathcal{D}_x^{(2)}$ hold for
every choice of $(\alpha_1,\alpha_2)$.

If $x$ is integer, we refer to Lemma \ref{lem:mgf-AbundoFuria-(3.17)-gen} and we have
$$\frac{1}{x}\log\mathbb{E}\left[e^{\alpha_1\tau(x)+\alpha_2\frac{A(x)}{x}}\right]
=\left\{\begin{array}{ll}
\frac{1}{x}\sum_{k=1}^x\varphi\left(\alpha_1+\alpha_2\frac{k}{x}\right)&\ \mbox{if}\ (\alpha_1,\alpha_2/x)\in\mathcal{D}_x^{(1)}\\
\infty&\ \mbox{otherwise}.
\end{array}\right.$$
Moreover, when we take the limit as $x\to\infty$, the inequalities that define the sets $\mathcal{D}_x^{(1)}$
lead to the set $\mathcal{D}$; in fact, if $\alpha_2\geq 0$, then we have
$\alpha_1+\frac{\alpha_2}{x}\cdot x\in\mathcal{D}(\varphi)$ for each fixed $x>0$, and therefore we get
$\alpha_1+\alpha_2\in\mathcal{D}(\varphi)$ because $x$ has no influence in the limit; if $\alpha_2<0$,
then we have $\alpha_1+\frac{\alpha_2}{x}\in\mathcal{D}(\varphi)$ for each fixed $x>0$, and therefore
we get $\alpha_1\leq\bar{\alpha}$ in the limit. In conclusion the limit coincides with the function
$\Lambda$ in Definition \ref{def:function-Lambda} because we trivially have the limit of an integral sum.

If $x$ is not integer, we refer to Lemma \ref{lem:mgf-AbundoFuria-(3.18)-gen} and we have
$$\frac{1}{x}\log\mathbb{E}\left[e^{\alpha_1\tau(x)+\alpha_2\frac{A(x)}{x}}\right]
=\left\{\begin{array}{ll}
\frac{1}{x}\sum_{k=0}^{[x]}\varphi\left(\alpha_1+\alpha_2\left(1-\frac{k}{x}\right)\right)
&\ \mbox{if}\ (\alpha_1,\alpha_2/x)\in\mathcal{D}_x^{(2)}\\
\infty&\ \mbox{otherwise}.
\end{array}\right.$$
Moreover, when we take the limit as $x\to\infty$, the inequalities that define the sets $\mathcal{D}_x^{(2)}$
lead to the set $\mathcal{D}$; in fact, if $\alpha_2\geq 0$, then we have
$\alpha_1+\frac{\alpha_2}{x}\cdot x\in\mathcal{D}(\varphi)$ for each fixed $x>0$, and therefore we get
$\alpha_1+\alpha_2\in\mathcal{D}(\varphi)$ because $x$ has no influence in the limit; if $\alpha_2<0$, then
we have $\alpha_1+\frac{\alpha_2}{x}(x-[x])\in\mathcal{D}(\varphi)$ for each fixed $x>0$, and therefore we get
$\alpha_1\leq\bar{\alpha}$ in the limit. In conclusion the limit coincides with the function
$\Lambda$ in Definition \ref{def:function-Lambda} because we have the limit of an integral sum; in
fact we have to consider
$$\frac{[x]+1}{x}\cdot\frac{1}{[x]+1}\sum_{k=0}^{[x]}\varphi\left(\alpha_1+\alpha_2\left(1-\frac{k}{x}\right)\right),$$
where $\frac{[x]+1}{x}\to 1$ and, by noting that $\frac{k}{[x]+1}\leq\frac{k}{x}<\frac{k+1}{[x]+1}$
for $k\in\{0,1,\ldots,[x]\}$, the remaining part is an integral sum of $\int_0^1\varphi(\alpha_1+\alpha_2(1-w))dw$,
which coincides with $\int_0^1\varphi(\alpha_1+\alpha_2y)dy$ after the change of variable $y=1-w$.

Then the proposition is proved by an application of Theorem \ref{th:GE}; in fact $(0,0)$ trivially
belongs to the interior of $\mathcal{D}(\Lambda)$, and the function $\Lambda$ is differentiable
throughout the interior of $\mathcal{D}(\Lambda)$ by Lemma \ref{lem:differentiability-of-Lambda}.
\end{proof}

We present some remarks on Proposition \ref{prop:LD}. We start with a brief discussion on the almost sure convergence
of the involved random variables.

\begin{remark}\label{rem:convergences}
We have $\Lambda^*(z_1,z_2)=0$ if and only if
$$(z_1,z_2)=\left(\frac{\partial\Lambda}{\partial\alpha_1}(0,0),\frac{\partial\Lambda}{\partial\alpha_2}(0,0)\right)
=\left(\varphi^\prime(0),\frac{1}{2}\varphi^\prime(0)\right).$$
(see eq. \eqref{eq:gradient-of-Lambda-at-the-origin} for the second equality). So, if we set
$$\Lambda^*(B^c):=\inf_{(z_1,z_2)\in B^c}\Lambda^*(z_1,z_2),$$
where $B$ is a sufficiently small open neighborhood of
$\left(\varphi^\prime(0),\frac{1}{2}\varphi^\prime(0)\right)$, we have
$\Lambda^*(B^c)>0$ and, for every $\eta\in(0,\Lambda^*(B^c))$,
$$P\left(\left(\frac{\tau(x)}{x},\frac{A(x)}{x^2}\right)\in B^c\right)\leq e^{-x(\Lambda^*(B^c)-\eta)}\ \mbox{for $x$ large enough}.$$
In conclusion $\left(\frac{\tau(x)}{x},\frac{A(x)}{x^2}\right)$ converges to
$\left(\varphi^\prime(0),\frac{1}{2}\varphi^\prime(0)\right)$ almost surely
by a standard application of Borel Cantelli Lemma.

We also remark that, if $x$ is integer, the almost sure convergence of 
$\frac{\tau(x)}{x}$ and $\frac{A(x)}{x^2}$ as $x\to\infty$ can be seen as a consequence of the law of 
large numbers; indeed we have
$$\frac{\tau(x)}{x}=\frac{T_1+\cdots+T_x}{x}\to\mathbb{E}[T_1]=\varphi^\prime(0),$$
which yields
$$\frac{A(x)}{x^2}=\frac{\sum_{k=1}^xT_k(x+k-1)}{x^2}\to\frac{\mathbb{E}[T_1]}{2}=\frac{\varphi^\prime(0)}{2}.$$
\end{remark}

Now we present some comments on the LDPs for the marginal random variables and on the expression 
of the rate function $\Lambda^*$ in Proposition \ref{prop:LD}.

\begin{remark}\label{rem:marginalLDPs-Chaganty}
Assume that Proposition \ref{prop:LD} provides a full LDP. Then we can obtain
the full LDPs of first and second components separately by standard applications
of the contraction principle (see e.g. Theorem 4.2.1 in \cite{DemboZeitouni}),
with rate functions $I_1$ and $I_2$ defined by
\begin{equation}\label{eq:marginal-rf}
I_1(z_1):=\inf_{z_2\in\mathbb{R}}\Lambda^*(z_1,z_2)\ \mbox{and}\ I_2(z_2):=\inf_{z_1\in\mathbb{R}}\Lambda^*(z_1,z_2).
\end{equation}
Moreover, if we set $J(z_2|z_1):=\Lambda^*(z_1,z_2)-I_1(z_1)$ (and we have
$J(z_2|z_1)\geq 0$ by the first equality in eq. \eqref{eq:marginal-rf}), then
\begin{equation}\label{eq:Chaganty-expression}
\Lambda^*(z_1,z_2)=J(z_2|z_1)+I_1(z_1).
\end{equation}
This equality has some analogies with the formula in Theorem 2.3 in
\cite{Chaganty} concerning large deviations for joint distributions; in such
a case $J(\cdot|z_1)$ can be interpreted as the rate function for the
conditional distributions of the second component given the first one.
In connection with this fact, in Section \ref{sec:ex} we consider the case of
Poisson process (i.e. the case where the function $\varphi$ is defined by eq.
\eqref{eq:varphi-exponential-case}), and we present an alternative proof of
Proposition \ref{prop:LD} (when $x$ is integer) based on the application of
Theorem 2.3 in \cite{Chaganty}.

Finally an application of Theorem \ref{th:GE} for $\left\{\frac{\tau(x)}{x}:x>0\right\}$
yields the following alternative expression of $I_1$:
\begin{equation}\label{eq:def-varphi-star}
I_1(z_1)=\sup_{\alpha_1\in\mathbb{R}}\{\alpha_1z_1-\Lambda(\alpha_1,0)\}
=\sup_{\alpha_1\in\mathbb{R}}\{\alpha_1z_1-\varphi(\alpha_1)\}=:\varphi^*(z_1).
\end{equation}
This is not surprising if we consider the LDP of $\left\{\frac{\tau(x)}{x}:x>0\right\}$
with $x$ integer; in fact a standard application of Cram\'er Theorem on $\mathbb{R}$
(see e.g. Theorem 2.2.3 in \cite{DemboZeitouni}) provides the LDP with rate function
$\varphi^*$ in eq. \eqref{eq:def-varphi-star}.
\end{remark}

We conclude with some minor remarks: the rate function $\Lambda^*$ in Proposition \ref{prop:LD} is equal to
infinity outside of the first orthant $T$ (Remark \ref{rem:primo-ottante}) and a discussion on the possibility
to get a full LDP in Proposition \ref{prop:LD} (Remark \ref{rem:convex-analysis-argument}).

\begin{remark}\label{rem:primo-ottante}
We have $0\leq A(x)\leq x\tau(x)$ almost surely by construction; therefore
$$P\left(\left(\frac{\tau(x)}{x},\frac{A(x)}{x^2}\right)\in T\right)=1,\ \mbox{where}\ T:=\{(z_1,z_2)\in\mathbb{R}^2:0\leq z_2\leq z_1\}.$$
Then, if Proposition \ref{prop:LD} provides a full LDP,	the lower bound for the open
set $T^c$ yields
$$-\infty=\liminf_{x\to\infty}\frac{1}{x}\log P\left(\left(\frac{\tau(x)}{x},\frac{A(x)}{x^2}\right)\in T^c\right)
\geq-\inf_{(z_1,z_2)\in T^c}\Lambda^*(z_1,z_2);$$
so we conclude that $\Lambda^*(z_1,z_2)=\infty$ for $(z_1,z_2)\in T^c$.
\end{remark}

\begin{remark}\label{rem:convex-analysis-argument}
It is interesting to learn when Proposition \ref{prop:LD} provides a full 
LDP, i.e. when we can neglect the intersection with the exposed points in the lower bound for
opens sets. We already know that, by statement (c) in Theorem \ref{th:GE}, this happens
if the function $\Lambda$ is lower semi-continuous and essentially smooth; so, for
instance, the full LDP holds if $\mathcal{D}(\varphi)=\mathbb{R}$.
In some cases the function $\Lambda$ is essentially smooth and not lower semi-continuous;
for instance this happens in the case of Poisson process, i.e. the case where the
function $\varphi$ is defined by eq. \eqref{eq:varphi-exponential-case}, studied in Section
\ref{sec:ex}. However the full LDP holds if the image of $\nabla\Lambda$ is the interior
of the set $T$ in Remark \ref{rem:primo-ottante}; in fact, in such a case, the function
$\Lambda^*$ is strictly convex on each convex subset $C$ of the image of $\nabla\Lambda$
by Theorem 4.1.2 in \cite{HiriarturrutyLemarechal} (which can be stated even if the
function $f$ in that theorem is not lower semi-continuous).
\end{remark}

\section{Moderate deviations of $\left\{\left(\frac{\tau(x)}{x},\frac{A(x)}{x^2}\right):x>0\right\}$}\label{sec:MD}
In this section we study moderate deviations. This terminology is used for a class of LDPs
(see Proposition \ref{prop:MD}) where the random variables and the speed function depend
on some positive scaling factors $\{a_x:x>0\}$ (such that \eqref{eq:MD-conditions} holds),
and all these LDPs are governed by the same quadratic rate function $\Psi_\Lambda^*$ that
uniquely vanishes at the origin $(z_1,z_2)=(0,0)$ (see eq. \eqref{eq:Psi-Lambda-star-esplicita}
below). The involved random variables are
$\left\{\left(\frac{\tau(x)}{x},\frac{A(x)}{x^2}\right):x>0\right\}$
in Proposition \ref{prop:LD} with the centering terms
$\left(\frac{\partial\Lambda}{\partial\alpha_1}(0,0),\frac{\partial\Lambda}{\partial\alpha_2}(0,0)\right)$
in eq. \eqref{eq:gradient-of-Lambda-at-the-origin}, and multiplied by
the divergent scalar factor $\sqrt{xa_x}$. Now we are ready for the second main result in 
this paper. 
The proof consists of an application of Theorem \ref{th:GE} with $h=2$, and the function $f$ is 
a convex quadratic function $\Psi_\Lambda$ (see eq. \eqref{eq:def-Psi-Lambda}) which has some 
relationship with the function $\Lambda$ in Definition \ref{def:function-Lambda}.

\begin{proposition}\label{prop:MD}
For every family of positive numbers $\{a_x:x>0\}$ such that
\begin{equation}\label{eq:MD-conditions}
a_x\to 0\ \mbox{and}\ xa_x\to\infty
\end{equation}
holds, the family of random variables
$\left\{\sqrt{xa_x}\left(\frac{\tau(x)}{x}-\varphi^\prime(0),\frac{A(x)}{x^2}-\frac{1}{2}\varphi^\prime(0)\right):x>0\right\}$
satisfies the LDP with speed $1/a_x$ and good rate function $\Psi_\Lambda^*$
defined by
\begin{equation}\label{eq:Psi-Lambda-star-esplicita}
\Psi_\Lambda^*(z_1,z_2)=\frac{1}{2}\sum_{i,j=1}^2\hat{c}_{ij}z_iz_j
\end{equation}
where $C^{-1}=(\hat{c}_{ij})_{i,j\in\{1,2\}}$ is the matrix in eq.
\eqref{eq:inverse-of-hessian-matrix-of-Lambda-at-the-origin}.
\end{proposition}
\begin{proof}
We want to apply the G\"artner Ellis Theorem; so we have to show that
$$\lim_{x\to\infty}\underbrace{\frac{1}{1/a_x}\log\mathbb{E}\left[e^{\frac{\sqrt{xa_x}}{a_x}\left(\alpha_1\left(\frac{\tau(x)}{x}-\varphi^\prime(0)\right)
+\alpha_2\left(\frac{A(x)}{x^2}-\frac{1}{2}\varphi^\prime(0)\right)\right)}\right]}_{=:\Psi_\Lambda(x;\alpha_1,\alpha_2)}=\Psi_\Lambda(\alpha_1,\alpha_2)
\ (\mbox{for every}\ (\alpha_1,\alpha_2)\in\mathbb{R}^2),$$
where
\begin{equation}\label{eq:def-Psi-Lambda}
\Psi_\Lambda(\alpha_1,\alpha_2):=\frac{1}{2}\sum_{i,j=1}^2c_{ij}\alpha_i\alpha_j,
\end{equation}
and $C=(c_{ij})_{i,j\in\{1,2\}}$ is the matrix in eq.
\eqref{eq:hessian-matrix-of-Lambda-at-the-origin} (which depends on the function $\Lambda$ in Definition \ref{def:function-Lambda}).
In fact the function $\Psi_\Lambda$ is trivially essentially smooth and lower semi-continuous and, after some standard
computations (we omit the details), one can check that
$$\Psi_\Lambda^*(z_1,z_2):=\sup_{(\alpha_1,\alpha_2)\in\mathbb{R}^2}\{\alpha_1z_1+\alpha_2z_2-\Psi_\Lambda(\alpha_1,\alpha_2)\}$$
coincides with $\Psi_\Lambda^*(z_1,z_2)$ in eq. \eqref{eq:Psi-Lambda-star-esplicita}.

In what follows we take into account that
\begin{multline*}
\Psi_\Lambda(x;\alpha_1,\alpha_2):=
\frac{1}{1/a_x}\log\mathbb{E}\left[e^{\frac{\sqrt{xa_x}}{a_x}\left(\alpha_1\left(\frac{\tau(x)}{x}-\varphi^\prime(0)\right)
+\alpha_2\left(\frac{A(x)}{x^2}-\frac{1}{2}\varphi^\prime(0)\right)\right)}\right]\\
=a_x\left(\log\mathbb{E}\left[e^{\left(\alpha_1\frac{\tau(x)}{\sqrt{xa_x}}+\alpha_2\frac{A(x)}{x\sqrt{xa_x}}\right)}\right]
-\frac{1}{\sqrt{xa_x}}\left(\alpha_1x\varphi^\prime(0)+\alpha_2\frac{x}{2}\varphi^\prime(0)\right)\right).
\end{multline*}
Moreover, as in the proof of Proposition \ref{prop:LD}, we distinguish two cases:
$x$ integer, and $x$ not integer.

If $x$ is integer we refer to Lemma \ref{lem:mgf-AbundoFuria-(3.17)-gen}. Then we take $x$ large enough to have
$$(\alpha_1/\sqrt{xa_x},\alpha_2/(x\sqrt{xa_x}))\in\mathcal{D}_x^{(1)};$$
note that we can do this for every $(\alpha_1,\alpha_2)\in\mathbb{R}^2$ because
$xa_x\to\infty$. So, for those values of $x$, we have
$$\Psi_\Lambda(x;\alpha_1,\alpha_2)
=a_x\left(\sum_{k=1}^x\varphi\left(\frac{\alpha_1}{\sqrt{xa_x}}+\frac{\alpha_2}{x\sqrt{xa_x}}k\right)
-\frac{1}{\sqrt{xa_x}}\left(\alpha_1x\varphi^\prime(0)+\alpha_2\frac{x}{2}\varphi^\prime(0)\right)\right).$$
Moreover we take into account the Maclaurin formula of order 2 for the function $z\mapsto\varphi(z)$, i.e.
$$\varphi(z)=\varphi^\prime(0)z+\varphi^{\prime\prime}(0)\frac{z^2}{2}+\frac{R(z)}{6}z^3,
\ \mbox{where}\ R(z)=\varphi^{\prime\prime\prime}(\omega(z))\ \mbox{for}\ |\omega(z)|\in(0,|z|);$$
therefore, for $z$ close enough to zero, there exists $M>0$ such that $\frac{|R(z)|}{6}\leq M$.
Then we obtain
\begin{multline*}
\Psi_\Lambda(x;\alpha_1,\alpha_2)
=a_x\left(\varphi^\prime(0)\sum_{k=1}^x\left(\frac{\alpha_1}{\sqrt{xa_x}}+\frac{\alpha_2}{x\sqrt{xa_x}}k\right)
+\frac{\varphi^{\prime\prime}(0)}{2}\sum_{k=1}^x\left(\frac{\alpha_1}{\sqrt{xa_x}}+\frac{\alpha_2}{x\sqrt{xa_x}}k\right)^2\right.\\
\left.+\frac{1}{6}\sum_{k=1}^xR\left(\frac{\alpha_1}{\sqrt{xa_x}}+\frac{\alpha_2}{x\sqrt{xa_x}}k\right)
\left(\frac{\alpha_1}{\sqrt{xa_x}}+\frac{\alpha_2}{x\sqrt{xa_x}}k\right)^3
-\frac{1}{\sqrt{xa_x}}\left(\alpha_1x\varphi^\prime(0)+\alpha_2\frac{x}{2}\varphi^\prime(0)\right)\right),
\end{multline*}
which can be rearranged as follows
\begin{multline*}
\Psi_\Lambda(x;\alpha_1,\alpha_2)
=a_x\left(\varphi^\prime(0)\frac{\alpha_1x}{\sqrt{xa_x}}+\varphi^\prime(0)\frac{\alpha_2}{x\sqrt{xa_x}}\cdot\frac{x(x+1)}{2}\right.\\
+\frac{\varphi^{\prime\prime}(0)}{2}\left(\frac{\alpha_1^2}{a_x}+\frac{2\alpha_1\alpha_2}{x^2a_x}\cdot\frac{x(x+1)}{2}+
\frac{\alpha_2^2}{x^3a_x}\cdot\frac{x(x+1)(2x+1)}{6}\right)\\
\left.+\frac{1}{6}\sum_{k=1}^xR\left(\frac{\alpha_1}{\sqrt{xa_x}}+\frac{\alpha_2}{x\sqrt{xa_x}}k\right)
\left(\frac{\alpha_1}{\sqrt{xa_x}}+\frac{\alpha_2}{x\sqrt{xa_x}}k\right)^3
-\frac{1}{\sqrt{xa_x}}\left(\alpha_1x\varphi^\prime(0)+\alpha_2\frac{x}{2}\varphi^\prime(0)\right)\right);
\end{multline*}
thus
\begin{multline*}
\Psi_\Lambda(x;\alpha_1,\alpha_2)
=a_x\left(\varphi^\prime(0)\frac{\alpha_2}{2\sqrt{xa_x}}+\frac{\varphi^{\prime\prime}(0)}{2}
\left(\frac{\alpha_1^2}{a_x}+\frac{\alpha_1\alpha_2(x+1)}{xa_x}+
\frac{\alpha_2^2}{x^3a_x}\cdot\frac{x(x+1)(2x+1)}{6}\right)\right.\\
\left.+\frac{1}{6}\sum_{k=1}^xR\left(\frac{\alpha_1}{\sqrt{xa_x}}+\frac{\alpha_2}{x\sqrt{xa_x}}k\right)
\left(\frac{\alpha_1}{\sqrt{xa_x}}+\frac{\alpha_2}{x\sqrt{xa_x}}k\right)^3\right)\\
=\frac{\alpha_2a_x\varphi^\prime(0)}{2\sqrt{xa_x}}+\frac{\varphi^{\prime\prime}(0)}{2}
\left(\alpha_1^2+\frac{\alpha_1\alpha_2(x+1)}{x}+
\frac{\alpha_2^2}{x^3}\cdot\frac{x(x+1)(2x+1)}{6}\right)\\
+\frac{a_x}{6}\sum_{k=1}^xR\left(\frac{\alpha_1}{\sqrt{xa_x}}+\frac{\alpha_2}{x\sqrt{xa_x}}k\right)
\left(\frac{\alpha_1}{\sqrt{xa_x}}+\frac{\alpha_2}{x\sqrt{xa_x}}k\right)^3.
\end{multline*}
So we can say that
$$\lim_{x\to\infty}\Psi_\Lambda(x;\alpha_1,\alpha_2)
=\frac{\varphi^{\prime\prime}(0)}{2}\left(\alpha_1^2+\alpha_1\alpha_2+\frac{\alpha_2^2}{3}\right)=\Psi_\Lambda(\alpha_1,\alpha_2)$$
(see eq. \eqref{eq:def-Psi-Lambda} for the last equality); in fact the limit of the linear
term and of the quadratic terms can be easily checked while the last sum multiplied by
$\frac{a_x}{6}$ tends to zero because, for $x$ large enough, we have
\begin{multline*}
\left|\frac{a_x}{6}\sum_{k=1}^xR\left(\frac{\alpha_1}{\sqrt{xa_x}}+\frac{\alpha_2}{x\sqrt{xa_x}}k\right)
\left(\frac{\alpha_1}{\sqrt{xa_x}}+\frac{\alpha_2}{x\sqrt{xa_x}}k\right)^3\right|\\
\leq\frac{Ma_x}{(xa_x)^{3/2}}\sum_{k=1}^x\left|\alpha_1+\alpha_2\frac{k}{x}\right|^3
\leq\frac{Ma_x}{(xa_x)^{3/2}}\sum_{k=1}^x(|\alpha_1|+|\alpha_2|)^3
=\frac{M(|\alpha_1|+|\alpha_2|)^3}{(xa_x)^{1/2}}.
\end{multline*}

If $x$ is not integer we refer to Lemma \ref{lem:mgf-AbundoFuria-(3.18)-gen}. Then we take 
$x$ large enough to have
$$(\alpha_1/\sqrt{xa_x},\alpha_2/(x\sqrt{xa_x}))\in\mathcal{D}_x^{(2)};$$
note that we can do this for every $(\alpha_1,\alpha_2)\in\mathbb{R}^2$ because
$xa_x\to\infty$. So, for those values of $x$, we have
$$\Psi_\Lambda(x;\alpha_1,\alpha_2)
=a_x\left(\sum_{k=0}^{[x]}\varphi\left(\frac{\alpha_1}{\sqrt{xa_x}}+\frac{\alpha_2}{x\sqrt{xa_x}}(x-k)\right)
-\frac{1}{\sqrt{xa_x}}\left(\alpha_1x\varphi^\prime(0)+\alpha_2\frac{x}{2}\varphi^\prime(0)\right)\right).$$
Then we can repeat some computations presented above (where $x$ is integer) with some slight
modifications and we get
\begin{multline*}
\Psi_\Lambda(x;\alpha_1,\alpha_2)
=a_x\left(\frac{\alpha_1([x]+1-x)\varphi^\prime(0)}{\sqrt{xa_x}}
+\frac{\alpha_2\varphi^\prime(0)}{\sqrt{xa_x}}\left(\frac{[x]+1}{x}\left(x-\frac{[x]}{2}\right)-\frac{x}{2}\right)\right.\\
\left.+\frac{\varphi^{\prime\prime}(0)}{2}\left(\frac{\alpha_1^2([x]+1)}{xa_x}+\frac{2\alpha_1\alpha_2}{x^2a_x}\cdot([x]+1)\left(x-\frac{[x]}{2}\right)+
\frac{\alpha_2^2}{x^3a_x}\cdot\sum_{k=0}^{[x]}(x-k)^2\right)\right.\\
\left.+\frac{1}{6}\sum_{k=0}^{[x]}R\left(\frac{\alpha_1}{\sqrt{xa_x}}+\frac{\alpha_2}{x\sqrt{xa_x}}(x-k)\right)
\left(\frac{\alpha_1}{\sqrt{xa_x}}+\frac{\alpha_2}{x\sqrt{xa_x}}(x-k)\right)^3\right)\\
=\frac{\alpha_1a_x([x]+1-x)\varphi^\prime(0)}{\sqrt{xa_x}}
+\frac{\alpha_2a_x\varphi^\prime(0)}{\sqrt{xa_x}}\left(\frac{[x]+1}{x}\left(x-\frac{[x]}{2}\right)-\frac{x}{2}\right)\\
+\frac{\varphi^{\prime\prime}(0)}{2}\left(\frac{\alpha_1^2([x]+1)}{x}+\frac{2\alpha_1\alpha_2}{x^2}\cdot([x]+1)\left(x-\frac{[x]}{2}\right)+
\frac{\alpha_2^2}{x^3}\cdot\sum_{k=0}^{[x]}(x-k)^2\right)\\
+\frac{a_x}{6}\sum_{k=0}^{[x]}R\left(\frac{\alpha_1}{\sqrt{xa_x}}+\frac{\alpha_2}{x\sqrt{xa_x}}(x-k)\right)
\left(\frac{\alpha_1}{\sqrt{xa_x}}+\frac{\alpha_2}{x\sqrt{xa_x}}(x-k)\right)^3.
\end{multline*}
So we can say that
$$\lim_{x\to\infty}\Psi_\Lambda(x;\alpha_1,\alpha_2)
=\frac{\varphi^{\prime\prime}(0)}{2}\left(\alpha_1^2+\alpha_1\alpha_2+\frac{\alpha_2^2}{3}\right)=\Psi_\Lambda(\alpha_1,\alpha_2)$$
(see eq. \eqref{eq:def-Psi-Lambda} for the last equality); in fact the first linear term tends to zero because $[x]+1-x$
is bounded, the second linear term tends to zero noting that
\begin{multline*}
\frac{[x]+1}{x}\left(x-\frac{[x]}{2}\right)-\frac{x}{2}=\frac{([x]+1)(2x-[x])-x^2}{2x}\\
=\frac{2x[x]-[x]^2+2x-[x]-x^2}{2x}=-\frac{(x-[x])^2}{2x}+\frac{x-[x]}{2x}+\frac{1}{2}\to\frac{1}{2},
\end{multline*}
the limits of the quadratic terms can be easily computed noting that
$$\frac{1}{x^3}\cdot\sum_{k=0}^{[x]}(x-k)^2=
\frac{1}{x^3}\left(x^2([x]+1)-2x\frac{[x]([x]+1)}{2}+\frac{[x]([x]+1)(2[x]+1)}{6}\right)\to\frac{1}{3},$$
and the last sum multiplied by $\frac{a_x}{6}$ tends to zero because, for $x$ large enough, we have
\begin{multline*}
\left|\frac{a_x}{6}\sum_{k=0}^{[x]}R\left(\frac{\alpha_1}{\sqrt{xa_x}}+\frac{\alpha_2}{x\sqrt{xa_x}}(x-k)\right)
\left(\frac{\alpha_1}{\sqrt{xa_x}}+\frac{\alpha_2}{x\sqrt{xa_x}}(x-k)\right)^3\right|\\
\leq\frac{Ma_x}{(xa_x)^{3/2}}\sum_{k=0}^{[x]}\left|\alpha_1+\alpha_2\frac{x-k}{x}\right|^3
\leq\frac{Ma_x}{(xa_x)^{3/2}}\sum_{k=0}^{[x]}(|\alpha_1|+|\alpha_2|)^3
=\frac{M([x]+1)a_x(|\alpha_1|+|\alpha_2|)^3}{(xa_x)^{3/2}}.
\end{multline*}
\end{proof}

We present some remarks on Proposition \ref{prop:MD}. We start with some typical features on moderate 
deviations.

\begin{remark}\label{rem:MD-typical-feature}
Typically moderate deviations fill the gap between the two asymptotic
regimes. In the case of Proposition \ref{prop:MD} we mean what follows
(as $x\to\infty$):
\begin{itemize}
\item the convergence of $\left(\frac{\tau(x)}{x}-\varphi^\prime(0),\frac{A(x)}{x^2}-\frac{1}{2}\varphi^\prime(0)\right)$
to zero (which is equivalent to the convergence of
$\left(\frac{\tau(x)}{x},\frac{A(x)}{x^2}\right)$ to
$\left(\frac{\partial\Lambda}{\partial\alpha_1}(0,0),\frac{\partial\Lambda}{\partial\alpha_2}(0,0)\right)$
stated in Remark \ref{rem:convergences});
\item the weak convergence of $\sqrt{x}\left(\frac{\tau(x)}{x}-\varphi^\prime(0),\frac{A(x)}{x^2}-\frac{1}{2}\varphi^\prime(0)\right)$
to the centered Normal distribution with covariance matrix $C$ in eq.
\eqref{eq:hessian-matrix-of-Lambda-at-the-origin}.
\end{itemize}
Note that the first asymptotic regime concerns the case $a_x=\frac{1}{x}$,
while the second one concerns the case $a_x=1$; so, in both cases, one
condition in eq. \eqref{eq:MD-conditions} holds, and the other one fails.
\end{remark}

The asymptotic normality result stated in Remark \ref{rem:MD-typical-feature} allows to provide
two approximate confidence intervals for $\varphi^\prime(0)$ when $x$ is large. In fact, if we
denote the standard Normal distribution function by $\Phi(\cdot)$, we obtain the following
approximate confidence intervals at the level $\ell\in(0,1)$:
$$\frac{\tau(x)}{x}\pm\frac{\sqrt{\varphi^{\prime\prime}(0)}}{\sqrt{x}}\Phi^{-1}((1+\ell)/2)$$
(which is also a consequence of the Central Limit Theorem, at least when $x$ is integer), and
$$2\left(\frac{A(x)}{x^2}\pm\frac{\sqrt{\varphi^{\prime\prime}(0)}}{\sqrt{3x}}\Phi^{-1}((1+\ell)/2)\right).$$
We remark that the second interval is larger than the first one because $\frac{2}{\sqrt{3}}>1$.

We conclude with some minor remarks: we discuss the possibility to present a slight modification 
of Proposition \ref{prop:MD} by considering different centering constants (Remark
\ref{rem:MD-same-result-with-different-centering-ee}) and we show how we can recover some limits
computed in \cite{AbundoFuria} which actually holds not only for Poisson processes (Remark
\ref{rem:MD-asymptotic-covariance-matrix}).

\begin{remark}\label{rem:MD-same-result-with-different-centering-ee}
The statement of Proposition \ref{prop:MD} still holds with the random
variables
\begin{equation}\label{eq:MD-rvs-with-different-centering}
\left\{\sqrt{xa_x}\left(\frac{\tau(x)}{x}-\frac{\mathbb{E}[\tau(x)]}{x},\frac{A(x)}{x^2}-\frac{\mathbb{E}[A(x)]}{x^2}\right):x>0\right\}
\end{equation}
in place of the random variables
$\left\{\sqrt{xa_x}\left(\frac{\tau(x)}{x}-\varphi^\prime(0),\frac{A(x)}{x^2}-\frac{1}{2}\varphi^\prime(0)\right):x>0\right\}$
(and, in particular, we can consider an alternative version of Remark
\ref{rem:MD-typical-feature} with appropriate changes). This can be proved
with slight changes of the proof of Proposition \ref{prop:MD} presented above.
However the result for the random variables in eq.
\eqref{eq:MD-rvs-with-different-centering} can be obtained by combining the
result in Proposition \ref{prop:MD} and Theorem 4.2.13 in \cite{DemboZeitouni}.
In fact, after some computations, we can check the exponential equivalence condition (see e.g.
Definition 4.2.10 in \cite{DemboZeitouni}), i.e.
\begin{multline*}
\lim_{x\to\infty}\frac{1}{1/a_x}\log P\left(\sqrt{xa_x}\left\|
\left(\frac{\tau(x)}{x}-\frac{\mathbb{E}[\tau(x)]}{x},\frac{A(x)}{x^2}-\frac{\mathbb{E}[A(x)]}{x^2}\right)\right.\right.\\
\left.\left.-\left(\frac{\tau(x)}{x}-\varphi^\prime(0),\frac{A(x)}{x^2}-\frac{1}{2}\varphi^\prime(0)\right)\right\|>\delta\right)
=-\infty\ (\mbox{for every}\ \delta>0),
\end{multline*}
where $\|\cdot\|$ is the Euclidean norm in $\mathbb{R}^2$.
\end{remark}

\begin{remark}\label{rem:MD-asymptotic-covariance-matrix}
We can check that the matrix $C$ in eq. \eqref{eq:hessian-matrix-of-Lambda-at-the-origin}
can be seen as an asymptotic covariance matrix. In particular we shall consider a
generalized version of some formulas in \cite{AbundoFuria} that concern the case of
exponentially distributed holding times; so here we have $\varphi^{\prime\prime}(0)$
in place of $\frac{1}{\lambda^2}$. We have the following limits as $x\to\infty$:
$$x\mathrm{Var}\left[\frac{\tau(x)}{x}\right]=\frac{\mathrm{Var}[\tau(x)]}{x}=
\left\{\begin{array}{ll}
\frac{x\varphi^{\prime\prime}(0)}{x}&\ \mbox{if}\ x\ \mbox{is integer}\\
\frac{([x]+1)\varphi^{\prime\prime}(0)}{x}&\ \mbox{if}\ x\ \mbox{is not integer}
\end{array}
\right.\to\varphi^{\prime\prime}(0)=c_{11}$$
(here we consider a generalized version of eq. (3.6) in \cite{AbundoFuria});
$$x\mathrm{Var}\left[\frac{A(x)}{x^2}\right]=\frac{\mathrm{Var}[A(x)]}{x^3}=
\frac{[x]+1}{12x^3}(12x(x-[x])+2[x](2[x]+1))\varphi^{\prime\prime}(0)\to
\frac{\varphi^{\prime\prime}(0)}{3}=c_{22}$$
(here we consider a generalized version of eqs. (3.8) and (3.9) in \cite{AbundoFuria}
for first and second moments of $A(x)$; however an explicit expression of the variance
appears as a factor in the denominator in eq. (3.16) in \cite{AbundoFuria});
$$x\mathrm{Cov}\left(\frac{\tau(x)}{x},\frac{A(x)}{x^2}\right)=\frac{\mathrm{Cov}(\tau(x),A(x))}{x^2}=
\left\{\begin{array}{ll}
\frac{x(x+1)\varphi^{\prime\prime}(0)}{2x^2}&\ \mbox{if}\ x\ \mbox{is integer}\\
\frac{([x]+1)(2x-[x])\varphi^{\prime\prime}(0)}{2x^2}&\ \mbox{if}\ x\ \mbox{is not integer}
\end{array}
\right.\to\frac{\varphi^{\prime\prime}(0)}{2}=c_{12}$$
(here we consider a generalized version of eqs. (3.12) and (3.15) in \cite{AbundoFuria}).

Finally, by taking into account the computations above, we have
$$\frac{\mathrm{Cov}(\tau(x),A(x))}{\sqrt{\mathrm{Var}[\tau(x)]\mathrm{Var}[A(x)]}}
=\frac{x\mathrm{Cov}\left(\frac{\tau(x)}{x},\frac{A(x)}{x^2}\right)}
{\sqrt{x\mathrm{Var}\left[\frac{\tau(x)}{x}\right]x\mathrm{Var}\left[\frac{A(x)}{x^2}\right]}}
\to\frac{c_{12}}{\sqrt{c_{11}c_{22}}}=\frac{\sqrt{3}}{2}\ \mbox{as}\ x\to\infty;$$
so the limit of the correlation coefficient computed in \cite{AbundoFuria} (see just after
eq. (3.16)) holds not only if the holding times are exponentially distributed.
\end{remark}

\section{On the case of Poisson process}\label{sec:ex}
Throughout this section we consider the case of Poisson process with intensity $\lambda>0$; so the random 
variables $\{T_n:n\geq 1\}$ in eq. \eqref{eq:renewal-process} are exponentially distributed and the function
$\varphi$ is as in eq. \eqref{eq:varphi-exponential-case}. Our aim is to discuss some aspects of the function
$\Lambda$ in this case, and in particular we refer to Remark \ref{rem:convex-analysis-argument}. Furthermore
we present an alternative proof of Proposition \ref{prop:LD} (when $x$ is integer) based on Theorem 2.3 in 
\cite{Chaganty}; in particular we get a slightly different expression of the rate function.

\subsection{On the function $\Lambda$ and Remark \ref{rem:convex-analysis-argument}}\label{sub:function-remark-PP}
We start by computing the function $\Lambda$ for $(\alpha_1,\alpha_2)\in\mathcal{D}$.
Firstly, if $\alpha_2=0$, we know that
\begin{equation}\label{eq:Lambda-Poisson-case-alpha2=0}
\Lambda(\alpha_1,0)=\varphi(\alpha_1)=\left\{\begin{array}{ll}
\log\frac{\lambda}{\lambda-\alpha_1}&\ \mbox{if}\ \alpha_1<\lambda\\
\infty&\ \mbox{otherwise}.
\end{array}\right.
\end{equation}
On the other hand, if $\alpha_2\neq 0$, we have
\begin{multline}
\Lambda(\alpha_1,\alpha_2)=\frac{1}{\alpha_2}\int_{\alpha_1}^{\alpha_1+\alpha_2}\varphi(y)dy
=\frac{1}{\alpha_2}\int_{\alpha_1}^{\alpha_1+\alpha_2}\log\frac{\lambda}{\lambda-y}dy\\
=\log\lambda+1+\frac{1}{\alpha_2}\{(\lambda-\alpha_1-\alpha_2)\log(\lambda-\alpha_1-\alpha_2)
-(\lambda-\alpha_1)\log(\lambda-\alpha_1)\}.\label{eq:Lambda-Poisson-case-alpha2-neq-0}
\end{multline}

The function $\Lambda$ is not lower semi-continuous because, as we said in Section
\ref{sub:function-Lambda}, the function $\varphi$ is integrable in a left neighborhood
of $\bar{\alpha}=\lambda$, i.e. $\int_0^\lambda\varphi(y)dy<\infty$. In fact the condition
in eq. \eqref{eq:lsc-condition} fails because, for $\alpha_2^{(0)}>0$, we have
$$\liminf_{(\alpha_1,\alpha_2)\to(\lambda-\alpha_2^{(0)},\alpha_2^{(0)})}\Lambda(\alpha_1,\alpha_2)
=\log\lambda+1-\log\alpha_2^{(0)}<\infty$$
and
$$\Lambda(\lambda-\alpha_2^{(0)},\alpha_2^{(0)})=\infty.$$

We recall (see Remark \ref{rem:convex-analysis-argument}) that the full LDP holds if the image
of $\nabla\Lambda$ coincides with the interior of the set $T$ in Remark \ref{rem:primo-ottante}.
In what follows we show that this happens; more precisely we mean that, for $z_1>z_2>0$, the system
$$\left\{\begin{array}{l}
z_1=\frac{\partial\Lambda}{\partial\alpha_1}(\alpha_1,\alpha_2)\\
z_2=\frac{\partial\Lambda}{\partial\alpha_1}(\alpha_1,\alpha_2)
\end{array}\right.$$
has a unique solution
$(\hat{\alpha}_1(z_1,z_2),\hat{\alpha}_2(z_1,z_2))$. We have two cases.
\begin{itemize}
\item If $\alpha_2\neq 0$, then
$$\left\{\begin{array}{l}
z_1=\frac{1}{\alpha_2}\log\frac{\lambda-\alpha_1}{\lambda-\alpha_1-\alpha_2}\\
z_2=\frac{1}{\alpha_2}((\lambda-\alpha_1)z_1-1),
\end{array}
\right.$$
or equivalently (after some manipulations starting from the second equality)
$$\left\{\begin{array}{l}
\log\frac{\alpha_2z_2+1}{\alpha_2(z_2-z_1)+1}=\alpha_2z_1\\
\lambda-\alpha_1=\frac{\alpha_2z_2+1}{z_1}.
\end{array}
\right.$$
\item If $\alpha_2=0$, then
$$\left\{\begin{array}{l}
z_1=\frac{1}{\lambda-\alpha_1}\\
z_2=\frac{1}{2(\lambda-\alpha_1)}.
\end{array}
\right.$$
\end{itemize}
We start with the case $\alpha_2\neq 0$; indeed the case $\alpha_2=0$ has to be
considered if and only if $2z_2=z_1$. The left hand side in the first equation,
i.e.
\begin{equation}\label{eq:def-function-g}
g(\alpha_2):=\log\frac{\alpha_2z_2+1}{\alpha_2(z_2-z_1)+1},
\end{equation}
is defined for $\alpha_2\in\left(-\frac{1}{z_2},\frac{1}{z_1-z_2}\right)$;
it is an increasing function because
$$g^\prime(\alpha_2)=\frac{z_1}{(\alpha_2z_2+1)(\alpha_2(z_2-z_1)+1)}>0$$
and we have $g(0)=0$ and $g^\prime(0)=z_1$. Moreover, by taking into
account its second derivative
$$g^{\prime\prime}(\alpha_2)=-z_1\frac{2\alpha_2(z_2-z_1)+2z_2-z_1}{(\alpha_2z_2+1)^2(\alpha_2(z_2-z_1)+1)^2},$$
we can say that $g(\alpha_2)$ is concave if $\alpha_2<\frac{2z_2-z_1}{2z_2(z_1-z_2)}$,
is convex if $\alpha_2>\frac{2z_2-z_1}{2z_2(z_1-z_2)}$. In conclusion one can realize
that, if $2z_2-z_1\neq 0$, then there exists $\alpha_2^*\neq 0$ such that $g(\alpha_2)=\alpha_2z_1$
if and only if $\alpha_2\in\{0,\alpha_2^*\}$. On the other hand, if $2z_2-z_1=0$, then
we have $g(\alpha_2)=\alpha_2z_1$ if and only if $\alpha_2=0$, and we set
$\alpha_2^*=0$. In conclusion the unique solution of the system is
$$(\hat{\alpha}_1(z_1,z_2),\hat{\alpha}_2(z_1,z_2))=\left(\lambda-\frac{\alpha_2^*z_2+1}{z_1},\alpha_2^*\right).$$

In the next Figure \ref{fig} we consider three different examples. As we know for each
example the abscissa $\alpha_2^*$ of the intersection between $g(\alpha_2)$ in eq.
\eqref{eq:def-function-g} and $h(\alpha_2)=\alpha_2z_1$ has the same sign of
$\gamma=\frac{2z_2-z_1}{2z_2(z_1-z_2)}$; we mean that we can have $\alpha_2^*,\gamma<0$
or $\alpha_2^*,\gamma>0$, or $\alpha_2^*=\gamma=0$ (see cases $(a)$, $(b)$ and $(c)$ in
Figure \ref{fig}, respectively).

\begin{figure}[H]
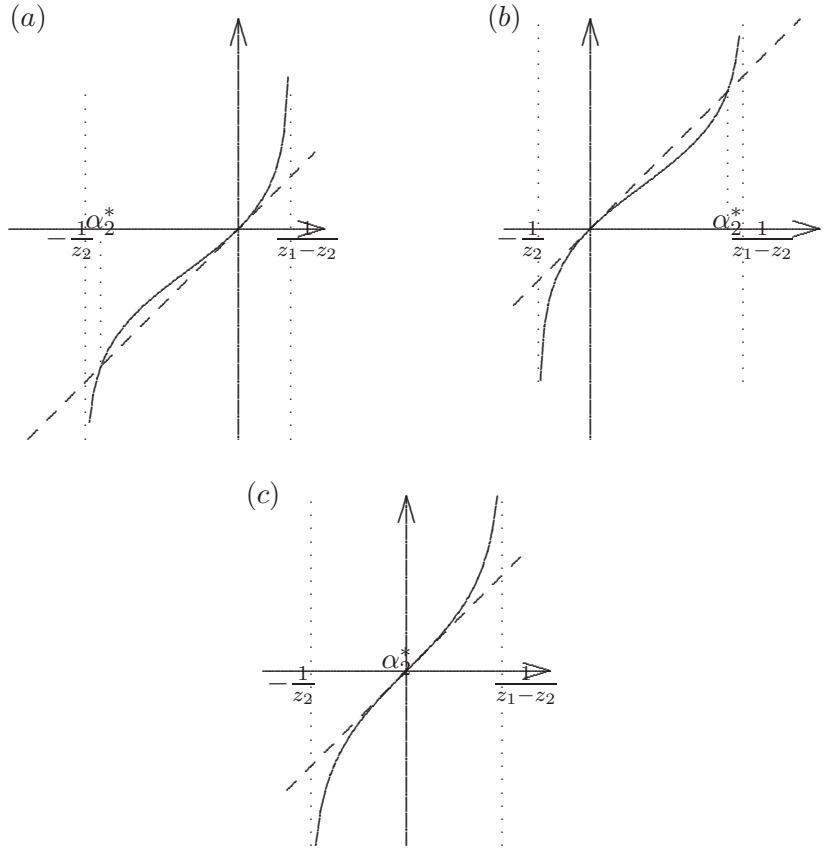

\hbox to\hsize\bgroup\hss

\beginpicture
\setcoordinatesystem units <0.2truein,0.2truein> \setplotarea x
from -3 to 3, y from -6 to 3

\plot
-3.9	-5.056
-3.8	-4.344
-3.7	-3.919
-3.6	-3.611
-3.5	-3.367
-3.4	-3.164
-3.3	-2.989
-3.2	-2.833
-3.1	-2.693
-3.0	-2.565
-2.9	-2.446
-2.8	-2.335
-2.7	-2.231
-2.6	-2.132
-2.5	-2.037
-2.4	-1.946
-2.3	-1.858
-2.2	-1.773
-2.1	-1.690
-2.0	-1.609
-1.9	-1.530
-1.8	-1.452
-1.7	-1.375
-1.6	-1.299
-1.5	-1.224
-1.4	-1.149
-1.3	-1.074
-1.2	-0.999
-1.1	-0.923
-1.0	-0.847
-0.9	-0.771
-0.8	-0.693
-0.7	-0.614
-0.6	-0.534
-0.5	-0.452
-0.4	-0.368
-0.3	-0.281
-0.2	-0.191
-0.1	-0.098
0.0	0.000
0.1	0.103
0.2	0.211
0.3	0.327
0.4	0.452
0.5	0.588
0.6	0.738
0.7	0.906
0.8	1.099
0.9	1.327
1.0	1.609
1.1	1.986
1.2	2.565
1.3 3.970  /

\put {$\alpha_2^*$} at -3.6 0.2
\put {$\frac 1{z_1-z_2}$} at 1.8 -0.3
\put {$-\frac 1{z_2}$} at -4.4 -0.3
\put {$(a)$} at -5.5 5.5
\arrow <10pt> [.3,.6] from 0 -5.5 to 0 5.5
\arrow <10pt> [.3,.6] from -6 0 to 2.25 0

\setdashes \plot -5.5 -5.5 2 2  /
\setdots   \plot -3.6 0  -3.6  -3.6  /
\plot 1.36  -5.5  1.36  3.5 /
\plot -4  -5.5  -4  3.5 /
\endpicture
\hspace{1.5cm}
\beginpicture
\setcoordinatesystem units <0.2truein,0.2truein> \setplotarea x
from -3 to 3, y from -6 to 6

\plot
-1.3	-3.970
-1.2	-2.565
-1.1	-1.986
-1.0	-1.609
-0.9	-1.327
-0.8	-1.099
-0.7	-0.906
-0.6	-0.738
-0.5	-0.588
-0.4	-0.452
-0.3	-0.327
-0.2	-0.211
-0.1	-0.103
0.0	0.000
0.1	0.098
0.2	0.191
0.3	0.281
0.4	0.368
0.5	0.452
0.6	0.534
0.7	0.614
0.8	0.693
0.9	0.771
1.0	0.847
1.1	0.923
1.2	0.999
1.3	1.074
1.4	1.149
1.5	1.224
1.6	1.299
1.7	1.375
1.8	1.452
1.9	1.530
2.0	1.609
2.1	1.690
2.2	1.773
2.3	1.858
2.4	1.946
2.5	2.037
2.6	2.132
2.7	2.231
2.8	2.335
2.9	2.446
3.0	2.565
3.1	2.693
3.2	2.833
3.3	2.989
3.4	3.164
3.5	3.367
3.6	3.611
3.7	3.919
3.8	4.344
3.9	5.056 /

\put {$\alpha_2^*$} at 3.6 0.2
\put {$\frac 1{z_1-z_2}$} at   4.5 -0.3
\put {$-\frac 1{z_2}$} at -1.8 -0.3
\put{$(b)$} at -2.25 5.5

\arrow <10pt> [.3,.6] from 0 -5.5  to 0 5.5
\arrow <10pt> [.3,.6] from -2.25   0 to 6 0

\setdashes \plot  -2 -2 5.5 5.5  /

\setdots \plot -1.36  -4 -1.36  5.55 /
\plot 4  -4 4  5.55 /
\plot 3.6 0  3.6  3.6 /
\endpicture

\hss\egroup \vglue8pt

\hbox to\hsize\bgroup\hss

\beginpicture
\setcoordinatesystem units <0.25truein,0.25truein> \setplotarea x
from 0 to 3, y from -3.5 to 3

\plot
-1.9	-3.664
-1.8	-2.944
-1.7	-2.512
-1.6	-2.197
-1.5	-1.946
-1.4	-1.735
-1.3	-1.551
-1.2	-1.386
-1.1	-1.237
-1.0	-1.099
-0.9	-0.969
-0.8	-0.847
-0.7	-0.731
-0.6	-0.619
-0.5	-0.511
-0.4	-0.405
-0.3	-0.302
-0.2	-0.201
-0.1	-0.100
0.0	0.000
0.1	0.100
0.2	0.201
0.3	0.302
0.4	0.405
0.5	0.511
0.6	0.619
0.7	0.731
0.8	0.847
0.9	0.969
1.0	1.099
1.1	1.237
1.2	1.386
1.3	1.551
1.4	1.735
1.5	1.946
1.6	2.197
1.7	2.512
1.8	2.944
1.9	3.664 /

\put {$\alpha_2^*$} at -0.2 0.2
\put {$\frac 1{z_1-z_2}$} at   2.5 -0.3
\put {$-\frac 1{z_2}$} at -2.4 -0.3
\put {$(c)$} at -3.0 3.66

\arrow <10pt> [.3,.6] from 0 -3.66 to 0 3.66
\arrow <10pt> [.3,.6] from -3 0 to 3 0

\setdashes \plot  -2.5 -2.5  2.5 2.5  /

\setdots  \plot -2  -3.66  -2  3.66 /
\plot 2  -3.66  2  3.66 /
\endpicture
\hss\egroup \vglue8pt

\caption{The functions $g$ (solid line) and $h$ (dashed line) for $z_1=1$ and three
different values of $z_2$. If $z_2=0.25$ (case $(a)$) then we have $2z_2<z_1$;
moreover $\left(-\frac{1}{z_2},\frac{1}{z_1-z_2}\right)=\left(-4,\frac{4}{3}\right)$
and $\alpha_2^*<0$. If $z_2=0.75$ (case $(b)$) then we have $2z_2>z_1$;
moreover $\left(-\frac{1}{z_2},\frac{1}{z_1-z_2}\right)=\left(-\frac{4}{3},4\right)$
and $\alpha_2^*>0$. If $z_2=0.5$ (case $(c)$) then we have $2z_2=z_1$;
moreover $\left(-\frac{1}{z_2},\frac{1}{z_1-z_2}\right)=(-2,2)$ and $\alpha_2^*=0$.}\label{fig}
\end{figure}

\subsection{The alternative proof of Proposition \ref{prop:LD} (when $x$ is integer)}\label{sub:alternative-proof}
We know that Proposition \ref{prop:LD} (together with the discussion above on Remark
\ref{rem:convex-analysis-argument}) yields the full LDP. Here we want to prove the
result with an application of Theorem 2.3 in \cite{Chaganty}. We recall that the
term \lq\lq proper rate function\rq\rq\ in \cite{Chaganty} coincides with the term
\lq\lq good rate function\rq\rq\ used in this paper. As we shall see we obtain a
different rate function expression $I_{\exp}$, say, in place of the rate function
$\Lambda^*$ in Proposition \ref{prop:LD} specified to the case of Poisson process,
i.e. the rate function defined by eq. \eqref{eq:def-Lambda-star}, where the function
$\Lambda$ is defined by eqs. \eqref{eq:Lambda-Poisson-case-alpha2=0} and
\eqref{eq:Lambda-Poisson-case-alpha2-neq-0}.

We shall consider an application of Theorem 2.3 in \cite{Chaganty} with
$\Omega_1=\Omega_2=[0,\infty)$. We have to check the following conditions:
\begin{enumerate}
\item the LDP of $\left\{\frac{\tau(x)}{x}:x>0\right\}$ holds with good rate function
$\varphi^*$, where $\varphi^*$ is defined by eq. \eqref{eq:def-varphi-star};
\item for a suitable family of good rate functions $\{\kappa^*(\cdot;z_1):z_1\geq 0\}$
(they will be presented in detail below), we have the LDP for the conditional distributions
$P\left(\frac{A(x)}{x^2}\in\cdot\Big|\frac{\tau(x)}{x}=z_1^{(x)}\right)$ as
$z_1^{(x)}\to z_1\in\Omega_1$ (as $x\to\infty$), with good rate function
$\kappa^*(\cdot;z_1)$;
\item the function $I$ defined by
\begin{equation}\label{eq:Chaganty-expression-Poisson}
I_{\exp}(z_1,z_2):=\varphi^*(z_1)+\kappa^*(z_2;z_1)
\end{equation}
is a good rate function.
\end{enumerate}

These three conditions will be checked below. Before doing this we introduce the
family of good rate functions $\{\kappa^*(\cdot;z_1):z_1\geq 0\}$ defined as follows:
\begin{equation}\label{eq:def-kappa-and-its-conjugate}
\kappa^*(z_2;z_1):=\sup_{\beta\in\mathbb{R}}\{\beta z_2-\kappa(\beta;z_1)\},\ \mbox{where}\
\kappa(\beta;z_1):=\left\{\begin{array}{ll}
\log\frac{e^{\beta z_1}-1}{\beta z_1}&\ \mbox{if}\ \beta\neq 0\\
0&\ \mbox{if}\ \beta=0.
\end{array}\right.
\end{equation}
These rate functions come up when one considers an application of Cram\'er
Theorem (already cited in the final part of Remark \ref{rem:marginalLDPs-Chaganty})
to obtain the LDP for the empirical means of i.i.d. random variables in $[0,z_1]$.
Note that, if $z_1=0$, we mean the trivial case of constant random variables
equal to zero, and therefore
$$\kappa^*(z_2;0):=\left\{\begin{array}{ll}
0&\ \mbox{if}\ z_2=0\\
\infty&\ \mbox{if}\ z_2\neq 0.
\end{array}\right.$$
On the contrary, if $z_1>0$ we do not have and explicit expression of
$\kappa^*(z_2;z_1)$; however we know that $\kappa^*(z_2;z_1)<\infty$ if $z_2\in(0,z_1)$
and $\kappa^*(z_2;z_1)=0$ if and only if $z_2=\frac{z_1}{2}$.

\paragraph{Condition 1.}
We already know (see the final part of Remark \ref{rem:marginalLDPs-Chaganty} where
we refer to Cram\'er Theorem on $\mathbb{R}$) that, when $x$ is integer,
$\left\{\frac{\tau(x)}{x}:x>0\right\}$ satisfies the LDP with rate function
$\varphi^*$. Actually it is easy to check that
$$\varphi^*(z_1)=\left\{\begin{array}{ll}
\lambda z_1-1-\log(\lambda z_1)&\ \mbox{if}\ z_1>0\\
\infty&\ \mbox{otherwise},
\end{array}\right.$$
that is a good rate function.

\paragraph{Condition 2.}
We want to apply  the G\"artner Ellis Theorem to the family of conditional distributions
of interest. Thus we have to check that
\begin{equation}\label{eq:GET-condition-for-conditional}
\lim_{x\to\infty}\frac{1}{x}\log\mathbb{E}\left[e^{x\beta\frac{A(x)}{x^2}}\Big|\frac{\tau(x)}{x}=z_1^{(x)}\right]
=\kappa(\beta;z_1)\ (\mbox{for every}\ \beta\in\mathbb{R}),\ \mbox{as}\ z_1^{(x)}\to z_1\in[0,\infty);
\end{equation}
actually the case $\beta=0$ can be neglected because it is trivial. We remark that
\begin{multline*}
\mathbb{E}\left[e^{\beta A(x)}\Big|\tau(x)=y\right]
=\int_0^\infty\cdots\int_0^\infty dt_1\cdots dt_ke^{\beta\sum_{k=1}^xkt_k}
\frac{\prod_{i=1}^k\lambda e^{-\lambda t_i}}{\frac{\lambda^x}{\Gamma(x)}y^{x-1}e^{-\lambda y}}1_{\{t_1+\cdots+t_k=y\}}\\
=\frac{\Gamma(x)}{y^{x-1}}\int_0^ydt_1\int_0^{y-t_1}dt_2\cdots\int_0^{y-(t_1+\cdots+t_{x-2})}dt_{x-1}e^{\beta\{\sum_{k=1}^{x-1}kt_k+x(y-(t_1+\cdots+t_{x-1}))\}}\\
=\frac{e^{\beta xy}(x-1)!}{y^{x-1}}
\underbrace{\int_0^ydt_1\int_0^{y-t_1}dt_2\cdots\int_0^{y-(t_1+\cdots+t_{x-2})}dt_{x-1}e^{-\beta\sum_{k=1}^{x-1}(x-k)t_k}}_{=:\mathcal{I}_x(\beta,y)},
\end{multline*}
where, for $\beta\neq 0$,
\begin{equation}\label{eq:formula-by-induction}
\mathcal{I}_x(\beta,y)=\frac{(1-e^{-\beta y})^{x-1}}{\beta^{x-1}(x-1)!}
\end{equation}
(see Appendix for details); therefore we obtain
$$\mathbb{E}\left[e^{\beta A(x)}\Big|\tau(x)=y\right]=\frac{e^{\beta xy}(x-1)!}{y^{x-1}}
\frac{(1-e^{-\beta y})^{x-1}}{\beta^{x-1}(x-1)!}=e^{\beta xy}\left(\frac{1-e^{-\beta y}}{\beta y}\right)^{x-1}.$$
Now we are ready to check the condition in eq. \eqref{eq:GET-condition-for-conditional}. As we said we neglect
the case $\beta=0$ and, since $z_1^{(x)}\to z_1\in[0,\infty)$ as $x\to\infty$, we have
\begin{multline*}
\frac{1}{x}\log\mathbb{E}\left[e^{x\beta\frac{A(x)}{x^2}}\Big|\frac{\tau(x)}{x}=z_1^{(x)}\right]
=\frac{1}{x}\log\left(e^{\beta xz_1^{(x)}}\left(\frac{1-e^{-\beta z_1^{(x)}}}{\beta z_1^{(x)}}\right)^{x-1}\right)\\
\to\beta z_1+\log\left(\frac{1-e^{-\beta z_1}}{\beta z_1}\right)
=\log\left(\frac{e^{\beta z_1}(1-e^{-\beta z_1})}{\beta z_1}\right)=\kappa(\beta;z_1)\ \mbox{as}\ x\to\infty.
\end{multline*}

\paragraph{Condition 3.}
Here we refer to Lemma 2.6 in \cite{Chaganty}, and we prove the goodness of the rate function
$I_{\exp}$ if we check the two following conditions.
\begin{itemize}
\item The function $(z_1,z_2)\mapsto\kappa^*(z_2;z_1)$ is lower semi-continuous.
\item For every compact subset $K_1$ of $\Omega_1$ and for every $L\geq 0$, the set
$$U(K_1,L):=\bigcup_{z_1\in K_1}\{z_2\in\Omega_2:\kappa^*(z_2;z_1)\leq L\}$$
is compact.
\end{itemize}
For the first condition we take $(z_1^{(n)},z_2^{(n)})\to (z_1,z_2)$ (as $n\to\infty$), and
we have
$$\kappa^*(z_2^{(n)};z_1^{(n)})\geq\beta z_2^{(n)}-\kappa(\beta;z_1^{(n)})$$
for all $\beta\in\mathbb{R}$, which yields the desired condition
$$\liminf_{n\to\infty}\kappa^*(z_2^{(n)};z_1^{(n)})\geq\kappa^*(z_2;z_1)$$
letting $n$ go to infinity (in fact $z_1\mapsto\kappa(\beta;z_1)$ is a continuous function)
and by taking the supremum with respect to $\beta$.

Now the second condition. We take a sequence $\{z_2^{(n)}:n\geq 1\}\subset U(K_1,L)$, and
therefore there exists a sequence $\{z_1^{(n)}:n\geq 1\}\subset K_1$ such that
$\kappa^*(z_2^{(n)};z_1^{(n)})\leq L$ for every $n\geq 1$. Then, by the compactness of
$K_1$, we can find a subsequence $\{z_1^{(n_k)}:k\geq 1\}\subset K_1$ such that
$z_1^{(n_k)}\to z_1^*\in K_1$. We show that $U(K_1,L)$ is compact if $\{z_2^{(n_k)}:k\geq 1\}$
(or a subsequence of $\{z_2^{(n_k)}:k\geq 1\}$) converges to a point of $z_2^*\in U(K_1,L)$.
We have two cases.\\
1) If $z_1^*=0$, then $z_2^{(n_k)}\to 0$ (because $\kappa^*(z_2^{(n)};z_1^{(n)})\leq L$ yields
$z_2^{(n)}\in[0,z_1^{(n)}]$) and
$$z_2^*:=0\in\{z_2\in\Omega_2:\kappa^*(z_2;0)\leq L\}\subset U(K_1,L).$$
2) If $z_1^*>0$, then we have
$$L\geq\kappa^*(z_2^{(n_k)};z_1^{(n_k)})=\sup_{\beta\in\mathbb{R}}\left\{\beta z_2^{(n_k)}-\kappa(\beta;z_1^{(n_k)})\right\}
=\sup_{\beta\in\mathbb{R}}\left\{\beta\frac{z_1^{(n_k)}}{z_1^*}z_2^{(n_k)}\frac{z_1^*}{z_1^{(n_k)}}
-\kappa\left(\beta\frac{z_1^{(n_k)}}{z_1^*};z_1^*\right)\right\},$$
and therefore
$$z_2^{(n_k)}\frac{z_1^*}{z_1^{(n_k)}}\in\{z_2\in\Omega_2:\kappa^*(z_2;z_1^*)\leq L\}.$$
Thus $\left\{z_2^{(n_k)}\frac{z_1^*}{z_1^{(n_k)}}:k\geq 1\right\}\subset\{z_2\in\Omega_2:\kappa^*(z_2;z_1^*)\leq L\}$,
i.e. we have a sequence of points in a compact set, and we can find a subsequence which
converges to a point $z_2^*$, say, such that
$$z_2^*\in\{z_2\in\Omega_2:\kappa^*(z_2;z_1^*)\leq L\}\subset U(K_1,L).$$
We conclude noting that $z_2^{(n_k)}\to z_2^*$ because $z_1^{(n_k)}\to z_1^*$.\\

So Conditions 1, 2 and 3 are checked. In particular, as a consequence of the discussion above, we can obtain
an interesting equality in terms of variational formulas. In fact $I_{\exp}$ in eq.
\eqref{eq:Chaganty-expression-Poisson} and $\Lambda^*$ in eq. \eqref{eq:Chaganty-expression} coincide when
$\varphi$ is as in eq. \eqref{eq:varphi-exponential-case}. Then, by taking into account eq.
\eqref{eq:def-varphi-star}, we get the equality
$$\kappa^*(z_2;z_1)=J(z_2;z_1)\ (\mbox{for all}\ z_1,z_2\geq 0),$$
where $\kappa^*(z_2;z_1)$ is given by \eqref{eq:def-kappa-and-its-conjugate}, and
$$J(z_2;z_1)=\sup_{(\alpha_1,\alpha_2)\in\mathbb{R}^2}\{\alpha_1z_1+\alpha_2z_2-\Lambda(\alpha_1,\alpha_2)\}
-(\lambda z_1-1-\log(\lambda z_1)).$$

\section{Conclusions}\label{sec:conclusions}
In this paper we have considered a renewal process $\{N(t):t\geq 0\}$ defined by \eqref{eq:renewal-process},
where $\{T_n:n\geq 1\}$ be i.i.d. light-tailed distributed (positive) random variables. Then, for $x\in\mathbb{R}_+$, 
we have considered the bivariate random variable $(\tau(x),A(x))$ where $\tau(x)$ is the first-passage time of 
$\{x-N(t):t\geq 0\}$ to reach zero or a negative value, and $A(x):=\int_0^{\tau(x)}(x-N(t))dt$ is the corresponding 
first-passage (positive) area swept out by the process $\{x-N(t):t\geq 0\}$. The asymptotic behavior of $A(x)$, as 
$x\to\infty$, has been already studied in the literature (also with some different hypotheses).

Our paper seems to be the first contribution for the study of the asymptotic behaviour of $(\tau(x),A(x))$.
We have presented some applications of the G\"artner Ellis Theorem; in particular we have studied moderate deviations,
i.e. a class of large deviation principles that fills the gap between a convergence to a constant vector of $\mathbb{R}^2$ 
and a weak convergence to the bivariate centered Normal distribution with a suitable covariance matrix. We have also given
some more details and an alternative proof of Proposition \ref{prop:LD} (when $x$ is integer) for the case in which 
$\{N(t):t\geq 0\}$ is a Poisson process.

In a future work one could try to obtain \emph{exact} asymptotic results as in the ones in \cite{PerfilevWatchel} for the
marginal distributions of $A(x)$ only. For instance one could try to consider the approach in \cite{Book} instead of an 
application the G\"artner Ellis Theorem.

\section*{Appendix: the proof of eq. \eqref{eq:formula-by-induction}}
Here we prove eq. \eqref{eq:formula-by-induction} by induction. We start with the case $x=2$ and we have
$$\mathcal{I}_2(\beta,y)=\int_0^ydt_1e^{-\beta\sum_{k=1}^{2-1}(2-k)t_k}=\int_0^ydt_1e^{-\beta t_1}=\frac{1-e^{-\beta y}}{\beta}.$$
So we assume that eq. \eqref{eq:formula-by-induction} is true, and we want to check that
it is also true for $x+1$. We have
\begin{multline*}
\mathcal{I}_{x+1}(\beta,y)=\int_0^ydt_1\int_0^{y-t_1}dt_2\cdots\int_0^{y-(t_1+\cdots+t_{x-1})}dt_xe^{-\beta\sum_{k=1}^x(x+1-k)t_k}\\
=\int_0^ydt_1\int_0^{y-t_1}dt_2\cdots\int_0^{y-(t_1+\cdots+t_{x-1})}dt_xe^{-\beta xt_1}e^{-\beta\sum_{k=2}^x(x+1-k)t_k}\\
=\int_0^ydt_1e^{-\beta xt_1}\int_0^{y-t_1}dt_2\cdots\int_0^{y-t_1-(t_2+\cdots+t_{x-1})}dt_xe^{-\beta\sum_{k=1}^{x-1}(x-k)t_k}
=\int_0^ydt_1e^{-\beta xt_1}\mathcal{I}_x(\beta,y-t_1);
\end{multline*}
so, by induction and some other computations, we get
\begin{multline*}
\mathcal{I}_{x+1}(\beta,y)=\int_0^ydt_1e^{-\beta xt_1}\frac{(1-e^{-\beta(y-t_1)})^{x-1}}{\beta^{x-1}(x-1)!}
=\int_0^ydt_1e^{-\beta t_1}\frac{(e^{-\beta t_1}-e^{-\beta y})^{x-1}}{\beta^{x-1}(x-1)!}\\
=\frac{1}{(-\beta)\beta^{x-1}(x-1)!}\left[\frac{(e^{-\beta t_1}-e^{-\beta y})^x}{x}\right]_{t_1=0}^{t_1=y}
=\frac{(1-e^{-\beta y})^x}{\beta^xx!}.
\end{multline*}
So eq. \eqref{eq:formula-by-induction} with $x+1$ is checked.

\paragraph{Acknowledgements.} The authors would like to thank the referees for their useful comments. We also thank Carlo 
Sinestrari for some discussion on Theorem 4.1.2 in \cite{HiriarturrutyLemarechal}.

\end{document}